\documentclass[mathpazo]{jiasc}
\usepackage[T1]{fontenc}
\usepackage{graphicx}
\usepackage{booktabs}
\usepackage{float}
\usepackage{placeins}
\usepackage{microtype}
\usepackage[hidelinks]{hyperref}
\hypersetup{
  pdftitle={Overcoming Spectral Bias via Cross-Attention},
  pdfauthor={Xiaodong Feng, Tao Tang, Xiaoliang Wan, Tao Zhou}
}

\renewcommand{\arraystretch}{1.5}

\begin{document}
\title{Overcoming Spectral Bias via Cross-Attention}
\author[X.~Feng, T.~Tang, X.~Wan and T.~Zhou]{Xiaodong Feng\affil{1},
  Tao Tang\affil{1,2}, Xiaoliang Wan\affil{3}, and
  Tao Zhou\affil{4}\comma\corrauth}
\address{\affilnum{1}\ Institute of Advanced Study, Beijing Normal-Hong Kong Baptist University, Zhuhai 519087, China.\\
  \affilnum{2}\ School of Mathematics and Statistics, Guangzhou Nanfang College, Guangzhou 510970, China.\\
  \affilnum{3}\ Department of Mathematics and Center for Computation and Technology, Louisiana State University, Baton Rouge 70803, USA\\
  \affilnum{4}\ SKLMS, Academy of Mathematics and Systems Science, Chinese Academy of Sciences}
\emails{\texttt{xiaodongfeng@bnbu.edu.cn} (X.~Feng),
  \texttt{ttang@nfu.edu.cn} (T.~Tang),
  \texttt{xlwan@lsu.edu} (X.~Wan),
  \texttt{tzhou@lsec.cc.ac.cn} (T.~Zhou)}

\begin{abstract}
Neural networks often exhibit spectral bias, a phenomenon where high-frequency components may converge substantially more slowly than low-frequency ones.
        To address this imbalance, we propose a cross-attention-based network that combines a multiscale random Fourier feature bank with learnable amplitude weighting.
    The amplitude envelope modulates the Fourier features according to their frequencies, while the cross-attention residual structure provides an input-dependent mechanism to emphasize informative scales. As a result, the proposed design accelerates high-frequency convergence relative to comparison networks built on the same multiscale bank. Moreover, the feature bank can be incrementally enriched with dominant modes identified through discrete Fourier analysis of intermediate approximations, without modifying the hidden-layer architecture. We further extend the framework to PDE learning through a two-network representation in which one subnetwork approximates higher-frequency solution components and the other approximates lower-frequency components, with a learnable or analytically computed mixing coefficient balancing their contributions. Numerical experiments on function approximation, image reconstruction, and representative elliptic PDEs demonstrate the accuracy of the proposed method for problems involving high-frequency, multiscale, and discontinuous features.
\end{abstract}
\maketitle

\section{Introduction}

In recent years, deep neural networks (DNNs) have been widely applied across various fields, including computer vision, speech recognition, natural language processing, and scientific computing for partial differential equations (PDEs) \cite{lecun2015deep,vaswani2017attention,liu2021swin,han2022survey,devlin2019bert,radford2021learning,raissi2019physics,karniadakis2021physics,weinan2018deep}.
Despite these advances, deploying conventional DNNs in computational science and engineering remains challenging.
A key limitation is their difficulty in capturing high-frequency content, a phenomenon often referred to as spectral bias or the Frequency Principle~\cite{xu2019frequency,rahaman2019spectral,xu2025overview}. DNNs
typically fit low-frequency components rapidly and generalize well at these frequencies, but struggle to represent high-frequency or highly oscillatory features \cite{xu2019training}.
This poses a significant obstacle for high-frequency and multiscale problems commonly encountered in PDE-based modeling.

Recent work has addressed this issue by incorporating explicit spectral
structure into DNNs~\cite{xu2025understanding,liu2024mitigating}. Approaches include frequency-aware initializations~\cite{zhang2025fourier}, choices of network structure~\cite{fang2024addressing,zhang2023shallow} and activation functions~\cite{jagtap2020adaptive,hong2022activation,sitzmann2020implicit} based on spectral properties, as well as augmentations with Fourier-type priors. For example, multiscale random Fourier feature (RFF) mappings, motivated by Neural Tangent Theory
\cite{jacot2018neural},  have been incorporated into physics-informed neural networks
(PINNs \cite{raissi2019physics}) to accelerate the convergence of high-frequency components
\cite{wang2021eigenvector,wang2022and}. Random feature methods \cite{chen2022bridging,chen2023random,chen2024optimization} can achieve spectral-like accuracy via multiscale representations and adaptive loss reweighting, particularly for PDEs on complex geometries.
The Fourier Multi-Component and Multi-Layer Neural Network (FMMNN) \cite{zhang2025fourier,zhang2025structured} aligns Fourier-structured activations with network design and employs scaled first-layer initialization to improve high-frequency approximation and optimization.
Alternatively, frequency manipulation can be performed at the solution level.
In \cite{cai2020phase}, a frequency-shifting approach uses phase shifts to relocate high-frequency content into lower-frequency regimes
during training and then reconstructs the original spectrum.
Moreover, multiscale networks such as MscaleDNN \cite{liu2020multi,wang2020multi} combine frequency scaling with
scale-separated subnetworks to represent multiscale structures more effectively. However, these methods typically rely on pre-specified frequency bases or prescribed multiscale transforms, which may limit adaptivity to instance-dependent spectral demands.
Beyond these representation and architectural strategies, recent work has shown that the choice of optimization method can also influence frequency-wise learning dynamics in physics-informed and operator learning~\cite{khodakarami2026spectral}. Hybrid strategies~\cite{lippe2023pde,oommen2025integrating,molinaro2024generative} that blend learning with iterative solvers, as well as multiscale neural operators and generative (e.g., diffusion-based) modeling, have also been explored to mitigate spectral bias.

More recently, adaptive frequency strategies \cite{xiong2025high,guibas2021adaptive,pathak2022fourcastnet} have been explored to reduce reliance on fixed multiscale designs.
In particular, frequency-adaptive multiscale methods \cite{huang2025frequency,huang2025frequency2} use discrete Fourier analysis of intermediate approximations to identify dominant scales and then recalibrate the multiscale mapping, leading to improved accuracy for PDEs with oscillatory solutions.
Moreover, frequency-guided PINNs (FG-PINNs) \cite{zheng2025fg} have been proposed to incorporate high-frequency prior information from the PDE and to separate and accelerate the learning of low- and high-frequency components.
While promising, FG-PINNs appear to depend strongly on the availability of informative high-frequency content in the source term, initial data, or boundary data, which may restrict their applicability in settings where such information is limited or not readily accessible.

This work argues that overcoming spectral bias requires not only richer frequency dictionaries, but also network operations that select and weight Fourier features according to the input and the evolving approximation. To this end, we develop a cross-attention-based design that adaptively reweights a learnably scaled multiscale Fourier feature bank and extend the same mechanism to PDE settings. Unlike self-attention mechanisms commonly used for sequence modeling, our approach employs cross-attention to adaptively select and weight Fourier features, allowing the network to emphasize different frequency bands in different regions of the domain.

\paragraph{Our goals and contributions}
Despite recent frequency-enriched approaches, the selection and emphasis of relevant frequencies are still often largely prespecified. Our goal is to introduce an attention-based mechanism that enables input-dependent spectral selection within a multiscale Fourier representation. Our contributions are summarized as follows:
\begin{itemize}
  \item We introduce a cross-attention architecture that uses a multiscale RFF bank with a learnable frequency-dependent amplitude envelope to perform input-dependent selection and reweighting of Fourier features across multiple scales.
  \item We design an adaptive frequency enhancement strategy that identifies dominant frequencies through discrete Fourier analysis of intermediate approximations and appends the corresponding modes to the multiscale feature bank.
  \item We extend the framework to PDE learning through a two-network formulation that combines the proposed cross-attention architecture with a fully connected network and uses a learnable or analytically computed mixing coefficient to balance their contributions.
\end{itemize}

The remainder of this paper is organized as follows. Section~\ref{sec:rff-ca} presents the multiscale RFF bank, its input-dependent reweighting through cross-attention, and the adaptive frequency enhancement strategy. Section~\ref{sec:PDE} extends the framework to PDE learning through a two-network formulation that balances low- and high-frequency solution components under physical constraints. Section~\ref{sec:numerical_experiments} reports numerical experiments on high-frequency, multiscale, and discontinuous function approximation, image-as-function approximation, and representative elliptic PDEs. Finally, Section~\ref{sec:conclusion} concludes the paper and discusses directions for future work.

\section{Cross-Attention for High-Frequency Function Approximation}
\label{sec:rff-ca}
In this section, we develop a neural network that uses cross-attention to incorporate information from a multiscale random Fourier feature (RFF) bank. We first describe the network structure, followed by a discussion on how to enhance high-frequency learning.
\subsection{Cross-Attention Network with a Multiscale Bank}
\label{subsec:multiscale-bank}
We propose a neural network for approximating high-frequency and discontinuous functions. The network incorporates multiscale Fourier features through cross-attention. In this design, a learnable amplitude envelope controls the feature amplitudes, while cross-attention performs input-dependent reweighting of the Fourier features.
\subsubsection{Bank of Multiscale Random Fourier Features}
Let $x \in \mathbb{R}^{d_{\mathrm{in}}}$ denote the input vector.
We construct a multiscale Fourier feature map $\phi(x)$ using randomly sampled frequencies and a learnable amplitude envelope.
The base frequencies are sampled as
\begin{equation*}
  \omega_m \sim \mathcal{N}\bigl(\mathbf{0},\, \sigma_\omega^{-2} I_{d_{\mathrm{in}}}\bigr),
  \qquad m=1,\ldots,M_{\mathrm{base}}.
\end{equation*}
The parameter $\sigma_\omega>0$ determines the characteristic scale of the base frequencies, with smaller values producing a broader frequency distribution. The sampled frequencies form the rows of a matrix $\Omega_{\mathrm{base}} \in
    \mathbb{R}^{M_{\mathrm{base}}\times d_{\mathrm{in}}}$. The same set of $M_{\mathrm{base}}$ base frequencies is used at every scale. For $k=0,\ldots,K$, multiplication by $2^k$ gives the frequencies associated with scale $k$:
\begin{equation*}
  \widetilde{\omega}_{m,k} \;=\; 2^{k}\omega_m.
\end{equation*}
Collecting all scaled frequencies as rows yields the multiscale frequency bank
\begin{equation*}
  \overline{\Omega}
  =
  \bigl[
      \Omega_{\mathrm{base}}^{\top},\,
      2\Omega_{\mathrm{base}}^{\top},\,
      \ldots,\,
      2^K\Omega_{\mathrm{base}}^{\top}
      \bigr]^{\top}
  \in \mathbb{R}^{M \times d_{\mathrm{in}}},
  \qquad M = M_{\mathrm{base}}(K+1).
\end{equation*}
For each scaled frequency $\widetilde{\omega}_{m,k}$, we sample a random phase as $b_{m,k} \sim \mathrm{Uniform}(0,2\pi).$
The sampled frequencies and phases remain unchanged during training.
To control the relative amplitudes of features at different frequencies, we introduce a learnable amplitude envelope:
\begin{equation*}
  a_{m,k} \;=\; \exp\bigl(-\beta \|\widetilde{\omega}_{m,k}\|_2 \bigr),
  \qquad \beta \geq 0,
\end{equation*}
where a single trainable scalar $\beta$ controls the rate of amplitude decay with frequency magnitude. When $\beta=0$, all features have equal amplitudes. For $\beta>0$, the amplitude weights decrease exponentially with frequency magnitude.

For each scale $k$, the scaled frequencies and their associated phases and amplitude weights define the Fourier feature vector
\begin{equation*}
  \phi_k(x)
    =
    \left[
      a_{m,k}\cos\bigl(\widetilde{\omega}_{m,k}^{\top}x+b_{m,k}\bigr)
      \right]_{m=1}^{M_{\mathrm{base}}}
    \in\mathbb{R}^{M_{\mathrm{base}}}.
\end{equation*}
Collecting the feature vectors from all $K+1$ scales gives
\begin{equation}
  \phi(x)
  \;=\;
  \frac{1}{\sqrt{M}}
    \bigl[\phi_0(x)^{\top},\,\phi_1(x)^{\top},\,\ldots,\,\phi_K(x)^{\top}\bigr]^{\top}
  \in \mathbb{R}^{M},
  \qquad M=M_{\mathrm{base}}(K+1).
\end{equation}

\subsubsection{Residual Updates with Cross-Attention}
To incorporate the multiscale Fourier features at each layer, we first organize them into groups. Let $d_q$ denote the network width and assume that $M$ is divisible by $d_q$. We partition $\phi(x)$ into $N_{\mathrm{base}}$ groups of $d_q$ consecutive features and arrange these groups as the rows of a matrix
\begin{equation}
  H(x)
  \in \mathbb{R}^{N_{\mathrm{base}}\times d_q},
  \qquad
  N_{\mathrm{base}}=\frac{M}{d_q}.
\end{equation}
We represent the intermediate network states as row vectors. Given an input $x$, the initial network representation is defined by
\begin{equation}\label{eqn:initial_feature}
  Q^{(0)}(x) \;=\; \sigma\bigl(\psi(x)^{\top} W^{(0)} + b^{(0)}\bigr)
  \in \mathbb{R}^{1\times d_q},
\end{equation}
where $\psi(x)=\phi(x)$ or $\psi(x)=x$, and $\sigma$ is an activation function applied componentwise. These two choices construct the initial representation from the Fourier features or directly from the input coordinates, respectively. The matrix $W^{(0)}$ has size $M\times d_q$ or $d_{\mathrm{in}}\times d_q$, respectively, and $b^{(0)}\in\mathbb{R}^{1\times d_q}$. Starting from $Q^{(0)}(x)$, we apply $L$ residual blocks to obtain the intermediate representations $Q^{(1)}(x),\ldots,Q^{(L)}(x)$. Each block incorporates the multiscale Fourier features supplied by the same matrix $H(x)$ through cross-attention.

Given an intermediate representation $Q^{(l)}(x)$ and the grouped Fourier features $H(x)$, the query, key, and value projections are defined as
\begin{equation*}
  Q_l = Q^{(l)}(x) W_Q^{(l)}, \qquad
  K_l = H(x) W_K^{(l)}, \qquad
  V_l = H(x) W_V^{(l)},
\end{equation*}
where $W_Q^{(l)}, W_K^{(l)}, W_V^{(l)} \in \mathbb{R}^{d_q \times d_q}$
are trainable matrices. Thus, $Q_l\in\mathbb{R}^{1\times d_q}$ and $K_l,V_l\in\mathbb{R}^{N_{\mathrm{base}}\times d_q}$.
The cross-attention output is
\begin{equation}
  \mathrm{CA}\bigl(Q^{(l)}(x), H(x)\bigr)
  \;=\;
  \mathrm{softmax}\left(\frac{Q_l K_l^{\top}}{\sqrt{d_q}}\right)\,V_l.
  \label{eq:ca}
\end{equation}
Here, softmax is applied across the $N_{\mathrm{base}}$ feature groups. The attention output is therefore a weighted combination of the projected feature groups in $V_l$, with weights that depend on $x$ through both $Q^{(l)}(x)$ and $H(x)$.

At each layer, the cross-attention output is added to the current representation through a residual connection:
\begin{equation}
  \widetilde{Q}^{(l)}(x)
  \;=\;
  Q^{(l)}(x) \;+\; \mathrm{CA}\bigl(Q^{(l)}(x), H(x)\bigr),
  \qquad l=0,1,\ldots,L-1.
  \label{eq:ca-residual}
\end{equation}
A nonlinear transformation with a second residual connection then gives the representation at the next layer:
\begin{equation}
  Q^{(l+1)}(x)
  \;=\;
  \widetilde{Q}^{(l)}(x) \;+\;
  \sigma\left(\widetilde{Q}^{(l)}(x)\,W_{\mathrm{ff}}^{(l)} + b_{\mathrm{ff}}^{(l)}\right),
  \label{eq:ffn-residual}
\end{equation}
where $W_{\mathrm{ff}}^{(l)}\in\mathbb{R}^{d_q\times d_q}$ and $b_{\mathrm{ff}}^{(l)}\in\mathbb{R}^{1\times d_q}$ are trainable parameters. After $L$ blocks, the network approximation is defined by
\begin{equation*}
  u(x;\theta) = Q^{(L)}(x)\,W_{\mathrm{out}} + b_{\mathrm{out}} .
\end{equation*}
Here, $\theta$ collects the trainable parameters. For scalar outputs, $W_{\mathrm{out}}\in\mathbb{R}^{d_q\times1}$ and $b_{\mathrm{out}}\in\mathbb{R}$.
The choices $\psi(x)=x$ and $\psi(x)=\phi(x)$ define the NN-CA and RFF-CA variants, respectively. For comparison, the corresponding networks without cross-attention are denoted NN and RFF-NN.

\begin{figure}[!htbp]
  \centering
  \includegraphics[width=0.9\linewidth]{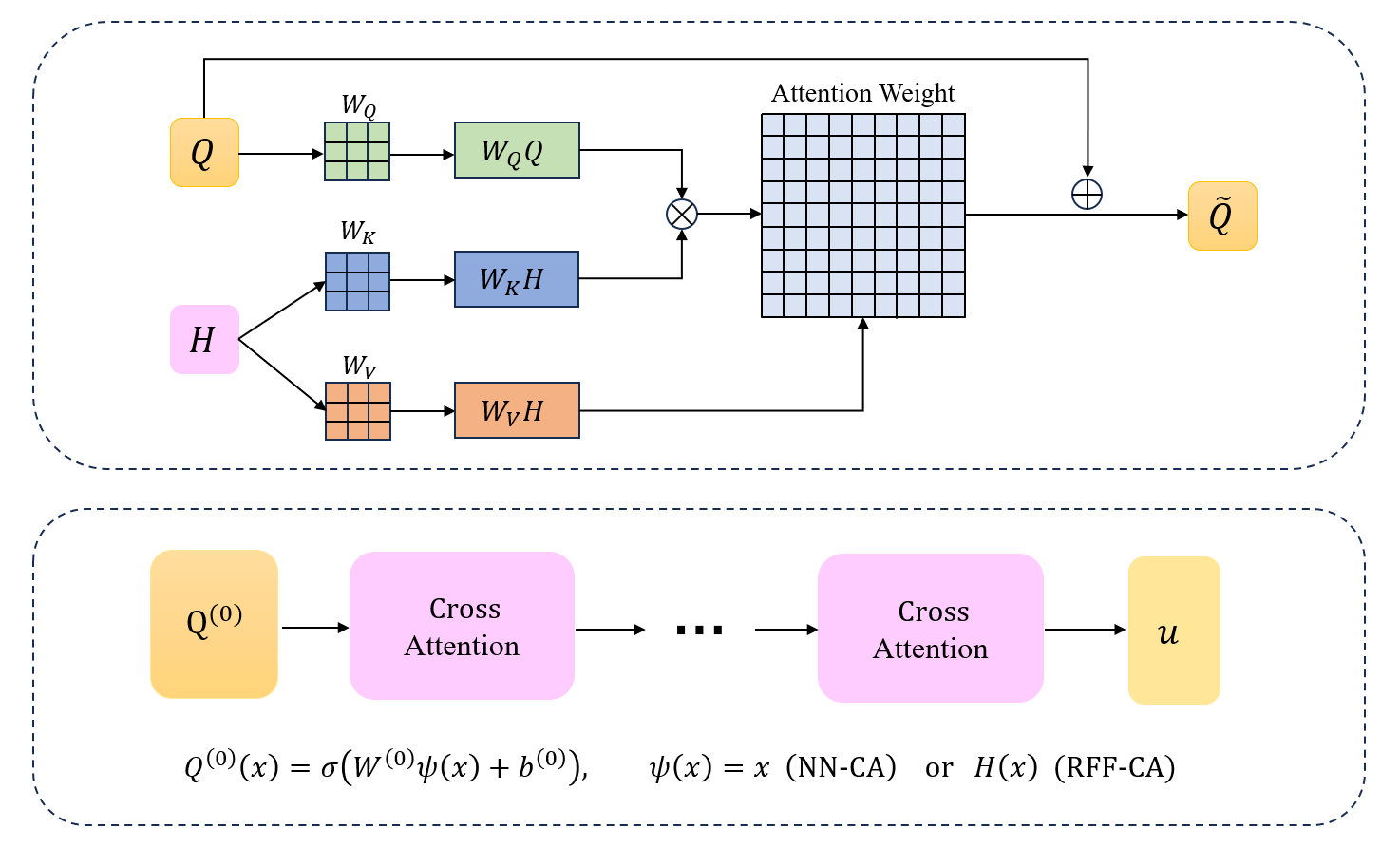}
  \caption{Schematic of the proposed NN-CA and RFF-CA network architectures.}
  \label{fig:network_schematic}
\end{figure}

From a spectral perspective, the amplitude envelope assigns frequency-dependent weights that are shared across the domain, while cross-attention provides input-dependent reweighting of the feature groups.
Figure~\ref{fig:network_schematic} illustrates the proposed network architecture.
\FloatBarrier
\subsection{Adaptive Frequency Enhancement}
\label{subsec:afe}

The construction in Section~\ref{subsec:multiscale-bank} provides input-dependent reweighting of multiscale Fourier features whose frequencies are chosen before training. To incorporate frequency information from a preliminary approximation, we introduce an adaptive frequency enhancement (AFE) strategy based on the discrete Fourier transform (DFT)~\cite{huang2025frequency}. The dominant frequencies identified from this approximation are used to construct additional Fourier features that enrich the initial bank. These features are incorporated into the cross-attention blocks using an additive attention mask.

Let $u(x;\theta^{(0)})$ denote a preliminary approximation obtained after training an NN-CA or RFF-CA model. We evaluate $u(x;\theta^{(0)})$ on a uniform grid over the unit periodic domain $\Omega=[0,1)^{d_{\mathrm{in}}}$ and apply a DFT to obtain the discrete Fourier coefficients $\widehat{u}_k(\theta^{(0)})$, indexed by a finite set $B\subset\mathbb{Z}^{d_{\mathrm{in}}}$. Other rectangular periodic domains can be treated by rescaling the coordinates.
Define
\begin{equation*}
  \zeta \;=\; \max_{k\in B}\bigl|\widehat{u}_k(\theta^{(0)})\bigr|,
\end{equation*}
and select the frequency indices whose coefficient magnitudes exceed a relative threshold:
\begin{equation}
  \mathcal{K}_{\mathrm{post}}
  \;=\;
  \Bigl\{\, k\in B:\ \bigl|\widehat{u}_k(\theta^{(0)})\bigr|>\lambda \zeta \Bigr\},
  \qquad 0<\lambda<1.
  \label{eq:afe-threshold}
\end{equation}
The parameter $\lambda$ controls the sparsity of the extracted modes.
We denote the original matrix of feature groups by
\begin{equation*}
  H_{\mathrm{base}}(x) \;=\; H(x)\in\mathbb{R}^{N_{\mathrm{base}}\times d_q}.
\end{equation*}
For each selected index $k\in\mathcal{K}_{\mathrm{post}}$, the corresponding frequency in unit-period coordinates is
\begin{equation*}
  \omega^{\mathrm{post}}_{k} = 2\pi k,
\end{equation*}
and these frequency vectors form the rows of a matrix $\Omega_{\mathrm{post}}\in\mathbb{R}^{M_{\mathrm{post}}\times d_{\mathrm{in}}}$, where $M_{\mathrm{post}}=|\mathcal{K}_{\mathrm{post}}|$. The additional Fourier feature vector is defined by
\begin{equation*}
  \phi_{\mathrm{post}}(x)
  \;=\;
  \sqrt{\frac{2}{M_{\mathrm{post}}}}\,
  \cos\bigl(\Omega_{\mathrm{post}}x + b_{\mathrm{post}}\bigr)
  \in\mathbb{R}^{M_{\mathrm{post}}},
\end{equation*}
where the cosine is applied componentwise and the entries of $b_{\mathrm{post}}\in\mathbb{R}^{M_{\mathrm{post}}}$ are sampled once from $\mathrm{Uniform}(0,2\pi)$.
We partition the entries of $\phi_{\mathrm{post}}(x)$ into groups of width $d_q$, padding the final group with zeros when necessary. These groups form the rows of a matrix
\begin{equation*}
  H_{\mathrm{post}}(x)\in\mathbb{R}^{N_{\mathrm{post}}\times d_q},
  \qquad
  N_{\mathrm{post}}=\left\lceil\frac{M_{\mathrm{post}}}{d_q}\right\rceil.
\end{equation*}
Stacking the original and additional feature groups by rows gives the augmented matrix
\begin{equation}
  H_{\mathrm{aug}}(x)
  \;=\;
  \begin{bmatrix}
      H_{\mathrm{base}}(x) \\
      H_{\mathrm{post}}(x)
    \end{bmatrix}
  \in\mathbb{R}^{(N_{\mathrm{base}}+N_{\mathrm{post}})\times d_q}.
  \label{eq:afe-H-aug}
\end{equation}
The augmented matrix is used to construct the keys and values in each cross-attention block:
\begin{equation*}
  K_l=H_{\mathrm{aug}}(x)W_K^{(l)},
    \qquad
    V_l=H_{\mathrm{aug}}(x)W_V^{(l)}.
\end{equation*}
The initial map $\psi(x)$ retains its definition in Section~\ref{subsec:multiscale-bank}, using the input coordinates or the original Fourier feature vector.

To control the contribution of the additional feature groups, we introduce an additive attention mask before applying softmax:
\begin{equation}
  A^{(l)}
  \;=\;
  \frac{Q_l K_l^{\top}}{\sqrt{d_q}}
  + \mathcal{M}^{(l)},
  \qquad
  \mathcal{M}^{(l)}
  =
  \bigl[
      \mathbf{0}_{1\times N_{\mathrm{base}}},\,
      \eta_l\,\mathbf{1}_{1\times N_{\mathrm{post}}}
      \bigr],
  \label{eq:afe-attn-mask}
\end{equation}
Both $A^{(l)}$ and $\mathcal{M}^{(l)}$ belong to $\mathbb{R}^{1\times(N_{\mathrm{base}}+N_{\mathrm{post}})}$. The zero block corresponds to the original feature groups, and the constant block corresponds to the additional groups. A negative value of $\eta_l$ suppresses the attention weights of the additional groups relative to those of the original groups. The resulting attention output is
\begin{equation}
  \mathrm{CA}\left(Q^{(l)}(x),H_{\mathrm{aug}}(x)\right)
  =
  \mathrm{softmax}\left(A^{(l)}\right)V_l
    \in\mathbb{R}^{1\times d_q}.
\end{equation}
We use a fixed value $\eta_l=\eta\leq 0$ in every layer throughout
  the second training stage and during evaluation.

In practice, AFE is applied as a two-stage refinement. A model is first trained with the original feature bank $H_{\mathrm{base}}(x)$. Frequency indices are then selected from the preliminary approximation using \eqref{eq:afe-threshold}, and training continues with the augmented feature bank $H_{\mathrm{aug}}(x)$ and the attention mask in \eqref{eq:afe-attn-mask}. The procedure may be repeated until $\mathcal{K}_{\mathrm{post}}$ stabilizes.
\begin{remark}[DFT on complex domains and in high dimensions]
  For a nonrectangular domain $\Omega$, one may extend the
    preliminary approximation to an enclosing box
    $\mathcal{Q}\supset\Omega$ and apply the DFT to uniform
    samples of the extended approximation.
    In high dimensions, when a tensor representation is
    available, one-dimensional DFTs of its component functions
    can be used to extract frequency information along each
    coordinate direction, avoiding a full multidimensional
    DFT \cite{huang2025frequency2}.
\end{remark}

\section{Application to numerical PDEs}
\label{sec:PDE}
\subsection{High-frequency amplification phenomenon}
\label{subsec:high-frequency-amplification}
As mentioned earlier, standard DNNs often exhibit a spectral bias toward learning low-frequency components, which can hinder the approximation of PDE solutions containing significant high-frequency content.
  When neural networks are trained with physics-based objectives, such as the PINN residual loss \cite{raissi2019physics} or the Deep Ritz energy \cite{weinan2018deep}, differential operators in the loss introduce a mode-dependent weighting of the approximation error.
In the Fourier domain, derivative terms act as wavenumber multipliers, thereby giving greater weight to high-frequency errors and potentially accelerating their decay during training. Related observations have also been reported in \cite{xu2025understanding,lu2021deepxde}.

To illustrate this effect, we consider the one-dimensional Poisson equation
\begin{equation*}
  -\Delta u(x)=f(x), \quad x\in\Omega=[-1,1],
\end{equation*}
with Dirichlet boundary conditions $u(-1)=u(1)=0$. We choose the exact solution
\begin{equation*}
  u(x)=\cos(\pi x)+\frac{1}{2}\cos(6\pi x)+\frac{1}{2}\cos(36\pi x).
\end{equation*}
We approximate $u(x)$ using an RFF-NN, consisting of a random Fourier feature map $\phi(x)$ followed by a fully connected network.
For each of the three objectives considered below, we use the same RFF map, network architecture, initial parameters, and $5000$ uniformly spaced training points.
Adam is used with a constant learning rate of $10^{-4}$; only the training objective is changed.

\begin{figure}[!htbp]
  \centering
  \includegraphics[width=\linewidth]{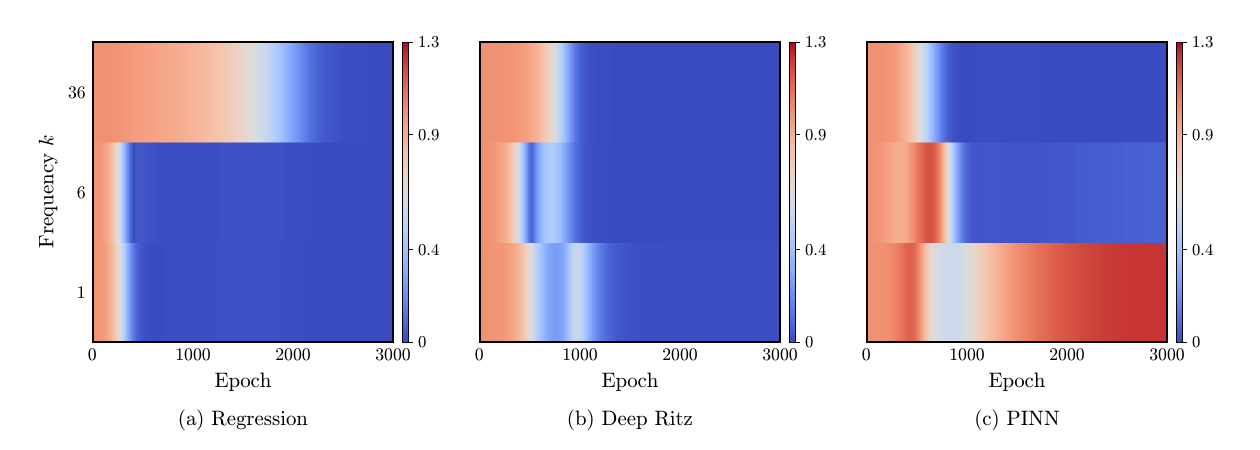}

  \caption{Training behavior of RFF-NN under the regression (left), Deep Ritz (middle), and PINN (right) objectives. Each heatmap reports the relative Fourier-coefficient error $\Delta_F(k)=\big|\widehat{u}^{\,\mathrm{pred}}_k-\widehat{u}_k\big|/|\widehat{u}_k|$ for the modes $k\in\{1,6,36\}$ over the first $3000$ epochs.}
  \label{fig:freq_principle}
\end{figure}

We compare three training objectives:
(i) a regression loss defined as the mean squared error (MSE) between the prediction and the exact solution
at the training points;
(ii) the Deep Ritz variational energy functional associated with the Poisson equation,
augmented by a penalty enforcing the Dirichlet boundary conditions;
and (iii) a PINN loss that penalizes the squared PDE residual $(-\Delta u(x;\theta)-f(x))$ together with the boundary mismatch.
Figure~\ref{fig:freq_principle} shows the evolution of the relative Fourier-coefficient errors at $k=1,6,36$ under the three objectives.

Under the regression loss, the relative Fourier-coefficient errors at $k=1$ and $k=6$ decrease earlier than that at $k=36$.
  Under the Deep Ritz objective, the error at $k=36$ decays on a time scale comparable to those of the lower modes.
  Under the PINN objective, the error at $k=36$ decreases before those of the lower modes, while the error at $k=1$ remains large over the displayed training interval.

This behavior is consistent with the Fourier weighting introduced by derivatives.
  For the basis function $\cos(k\pi x)$, the first-derivative term in the Deep Ritz objective contributes a modal weight proportional to $(k\pi)^2$, whereas the squared second-order residual in the PINN objective contributes a weight proportional to $(k\pi)^4$.
  These frequency-dependent weights can alter the relative decay rates of the modal errors during training.
  An idealized modal analysis of this effect is given in \ref{appendix:1}.

\subsection{Mitigating High-Frequency Amplification}
As illustrated in Section~\ref{subsec:high-frequency-amplification}, low-frequency errors may remain large even when high-frequency errors decay rapidly. To improve the approximation of low- and high-frequency components, we employ a linear combination of two subnetworks.
Specifically, consider a PDE with Dirichlet boundary conditions:
\begin{align*}
  \mathcal{N}[u](x) & = f(x), \quad x \in \Omega,          \\
  u(x)                        & = g(x), \quad x \in \partial \Omega.
\end{align*}
Here $\mathcal{N}[\cdot]$ is a differential operator, and $\Omega$ is the computational domain with boundary $\partial \Omega$. We approximate the solution by
\begin{equation*}
  u(x;\theta) = u_h(x;\theta_h) + \alpha \, u_l(x;\theta_l),
\end{equation*}
Here, $u_h$ is represented by an RFF-CA (or NN-CA) network to approximate oscillatory features of the solution. A fully connected network $u_l$ is introduced to improve the approximation of lower-frequency components, with the scalar coefficient $\alpha$ controlling its contribution to the combined approximation.
We impose the boundary conditions weakly through separate penalties on $u_h-g$ and $\alpha u_l$ on $\partial\Omega$. The training loss is then defined as
\begin{align*}
  \mathcal{L} = \int_{\Omega} \big(\mathcal{N}[u_h+\alpha \, u_l](x) - f(x)\big)^2 \rho_r(x) \mathrm{d}x + \gamma \int_{\partial \Omega} \big((u_h - g)^2 + \alpha^2 u_l^2\big) \rho_b(x) \, \mathrm{d}S,
\end{align*}
where $\rho_r$ and $\rho_b$ are sampling densities on $\Omega$ and $\partial\Omega$, respectively, and $\gamma>0$ weights the boundary penalties. We propose two strategies to determine the mixing coefficient $\alpha$:
\begin{itemize}
  \item The coefficient $\alpha$ is optimized jointly with the network parameters $\theta_h$ and $\theta_l$.
  \item For a linear differential operator $\mathcal{N}$, the loss $\mathcal{L}$ is quadratic in $\alpha$ for the current network parameters $\theta_h$ and $\theta_l$. At each training epoch, we compute its minimizer, which is unique provided that the denominator below is positive:
      \begin{align*}
      \alpha_{\mathrm{opt}} = -\frac{\int_{\Omega} \big(\mathcal{N}[u_h](x)-f(x)\big) \mathcal{N}[u_l](x)\rho_r(x)\, \mathrm{d}x}{\int_{\Omega} (\mathcal{N}[u_l])^2(x)\rho_r(x) \,\mathrm{d}x + \gamma \int_{\partial \Omega}u_l^2(x) \rho_b(x)\,\mathrm{d}S}.
      \end{align*}
      The integrals in $\mathcal{L}$ and in the expression for $\alpha_{\mathrm{opt}}$ are approximated by sample averages:
      \begin{equation*}
      \mathcal{L} \approx
      \frac{1}{N_r}\sum_{i=1}^{N_r}\big(\mathcal{N}[u_h+\alpha_{\mathrm{opt}} u_l]-f\big)^2\big(x_r^{(i)}\big)
      +
      \frac{\gamma}{N_b}\sum_{i=1}^{N_b}\Big( (u_h-g)^2+\alpha_{\mathrm{opt}}^{2}u_l^2 \Big)\big(x_b^{(i)}\big),
      \end{equation*}
      where
      \begin{equation*}
      \alpha_{\mathrm{opt}}\approx
      -\frac{\displaystyle \frac{1}{N_r}\sum_{i=1}^{N_r}(\mathcal{N}[u_h] -f)\big(x_r^{(i)}\big)\,\mathcal{N}[u_l]\big(x_r^{(i)}\big)}
      {\displaystyle \frac{1}{N_r}\sum_{i=1}^{N_r}\big(\mathcal{N}[u_l]\big)^2\big(x_r^{(i)}\big)
      +\frac{\gamma}{N_b}\sum_{i=1}^{N_b}\big(u_l\big)^2\big(x_b^{(i)}\big)} \, ,
      \end{equation*}
      and $\{x_r^{(i)}\}$ and $\{x_b^{(i)}\}$ are the collocation points sampled from $\rho_r$ and $\rho_b$, respectively.
      In practice, we remove the factor $\alpha_{\mathrm{opt}}^2$ from the boundary penalty on $u_l$, so that the strength of this penalty is independent of $\alpha_{\mathrm{opt}}$. The resulting loss for the network parameters is
      \begin{equation*}
      \widehat{\mathcal{L}} \approx
      \frac{1}{N_r}\sum_{i=1}^{N_r}\big(\mathcal{N}[u_h+\alpha_{\mathrm{opt}} u_l]-f\big)^2\big(x_r^{(i)}\big)
      +
      \frac{\gamma}{N_b}\sum_{i=1}^{N_b}\Big( (u_h-g)^2+u_l^2 \Big)\big(x_b^{(i)}\big).
      \end{equation*}

      The coefficient computed from $\mathcal{L}$ is treated as a constant when differentiating $\widehat{\mathcal{L}}$ with respect to $(\theta_h,\theta_l)$.

\end{itemize}

\begin{figure}[!htbp]
  \centering
  \includegraphics[width=0.88\linewidth]{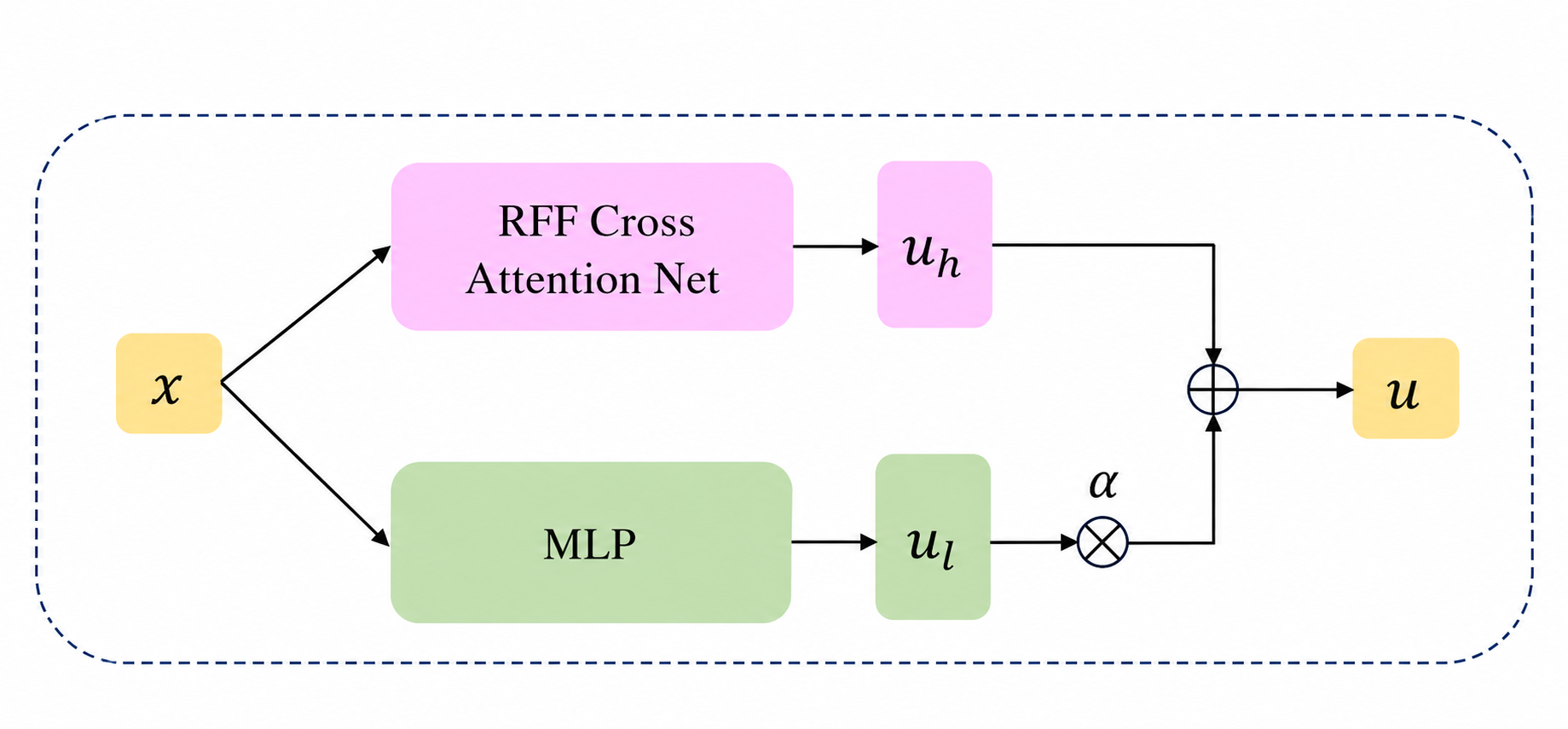}
  \caption{Schematic of the two-network PDE approximation using an RFF-CA network and a fully connected network.}
  \label{fig:pde_solver}
\end{figure}

A schematic of the proposed method is shown in Fig.~\ref{fig:pde_solver}.

\FloatBarrier
\section{Numerical experiments}
\label{sec:numerical_experiments}
We assess the proposed approximation methods through function regression and elliptic PDE problems. Unless otherwise stated, we use $M_{\mathrm{base}}=128$, $K=3$, network width $d_q=64$, four attention heads, and $L=4$ residual blocks, with the amplitude parameter initialized to $\beta_0=0.1$.
The default training configuration uses double precision, Adam, and gradient clipping at $1.0$.

For all experiments, we report the relative $L^2$ error, computed on the evaluation points specified in the corresponding subsections, as
\begin{equation*}
  \mathrm{Rel}\,L^2(u(\cdot;\theta),u)=\frac{\|u(\cdot;\theta)-u\|_2}{\|u\|_2},
\end{equation*}
where $u$ and $u(\cdot;\theta)$ denote the reference solution or function and its approximation, respectively.

\subsection{Function approximation}
We compare the approximation accuracy of RFF-NN, RFF-CA, and NN-CA for oscillatory and discontinuous functions and natural images. A periodic function with a sparse Fourier spectrum is used to assess the effect of adaptive frequency enhancement.
\subsubsection{High-frequency and discontinuous function approximation}
\begin{figure}[!htbp]
  \centering
  \includegraphics[width=0.82\linewidth]{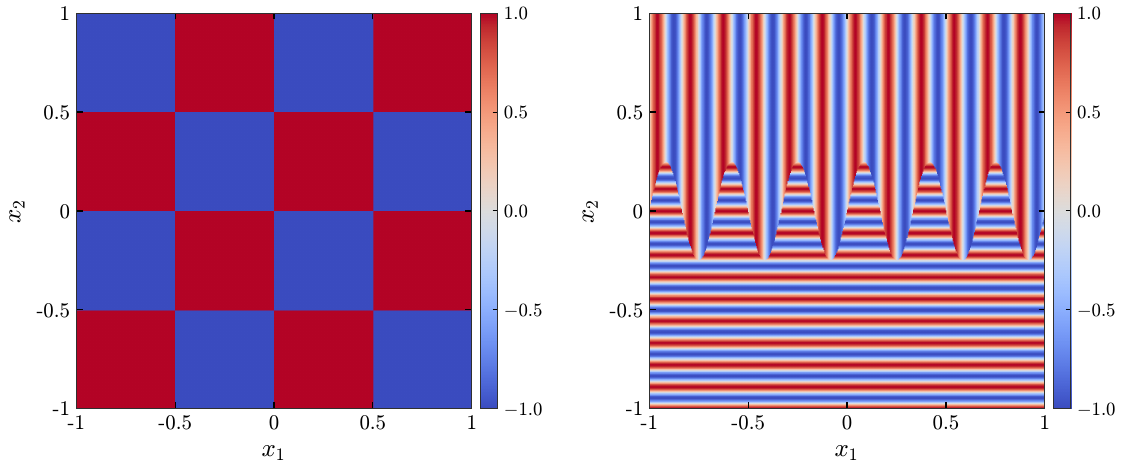}
  \caption{Functions $f_1$ (left) and $f_2$ (right).}
  \label{fig:regression_truth}
\end{figure}

We consider the following two functions on $\Omega=[-1,1]^2$.
  \begin{align*}
    f_1(x_1,x_2)
     & = \operatorname{sign}\bigl(\sin(2\pi x_1)\sin(2\pi x_2)\bigr), \\
    f_2(x_1,x_2)
     & = \begin{cases}
      \sin(12\pi x_1), & x_2>\tfrac14\sin(6\pi x_1),    \\
      \cos(18\pi x_2), & x_2\leq\tfrac14\sin(6\pi x_1).
      \end{cases}
  \end{align*}
As illustrated in Fig.~\ref{fig:regression_truth}, the checkerboard function $f_1$ has axis-aligned jumps, whereas $f_2$ combines oscillations in two directions across a curved interface.

We compare RFF-NN, RFF-CA, and NN-CA using the same initial Fourier feature bank within each example. The bank comprises $128$ base frequencies at one scale for $f_1$ and $256$ base frequencies at four dyadic scales for $f_2$, with $\sigma_\omega=1$ in both cases.

We train all models on a uniform $500\times500$ grid including the boundary of $\Omega$, using randomly shuffled mini-batches of size $4000$. Optimization uses Adam with no weight decay, starting from a learning rate of $2\times10^{-3}$ that is halved every $50$ epochs for $f_1$ and every $100$ epochs for $f_2$. The models are trained for $500$ and $1000$ epochs for $f_1$ and $f_2$, respectively, with the mean squared and relative $L^2$ errors evaluated on the entire training grid after each epoch.

\begin{figure}[!htbp]
  \centering
  \includegraphics[width=0.9\linewidth]{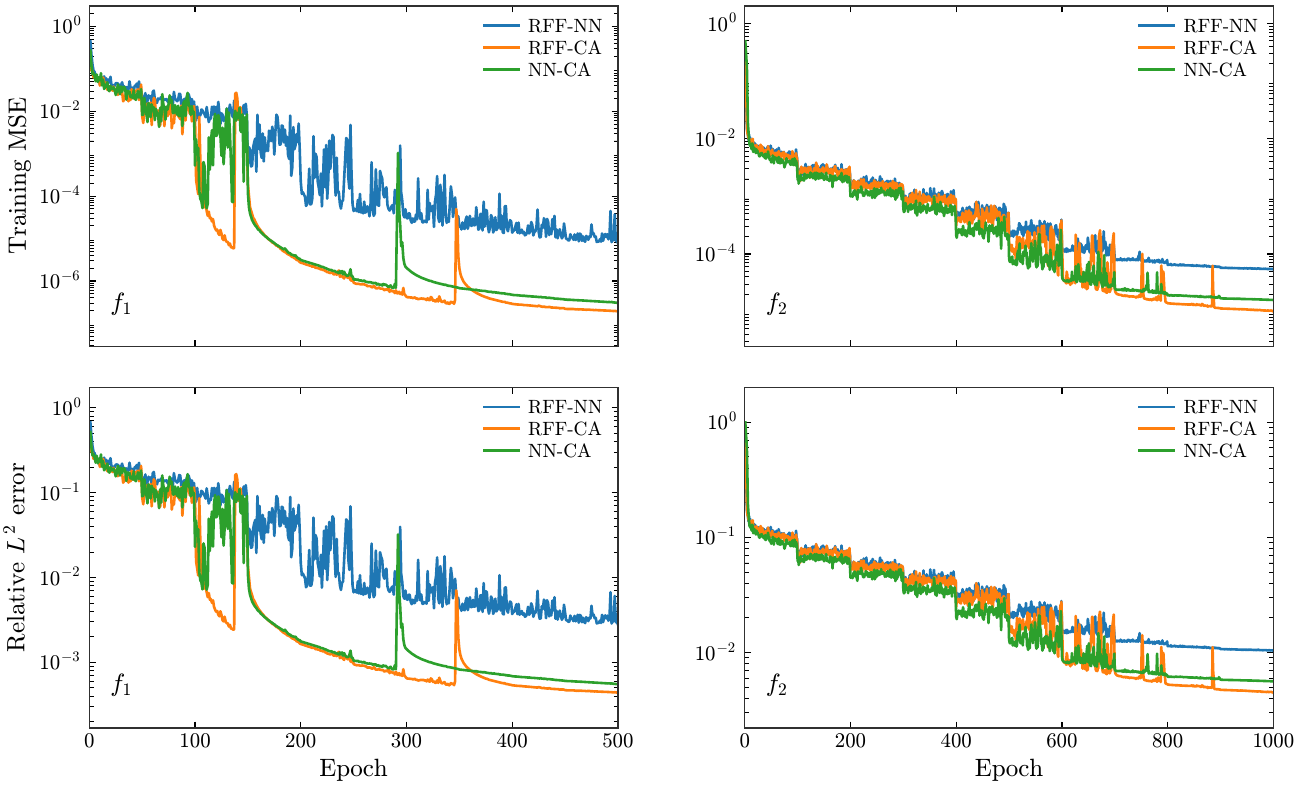}
  \caption{The mean squared error (top) and relative $L^2$ error (bottom) are shown during training for $f_1$ (left) and $f_2$ (right), with both errors evaluated on the training grid.}
  \label{fig:regression_results}
\end{figure}

Both RFF-CA and NN-CA attain smaller final training-grid errors than RFF-NN, as shown in Fig.~\ref{fig:regression_results}. RFF-CA gives the lowest relative $L^2$ error in both examples, reducing the error of RFF-NN by factors of approximately $6.4$ and $2.3$ for $f_1$ and $f_2$, respectively (Table~\ref{tab:synthetic_regression_errors}).

\begin{figure}[!htbp]
  \centering
  \includegraphics[width=0.95\linewidth]{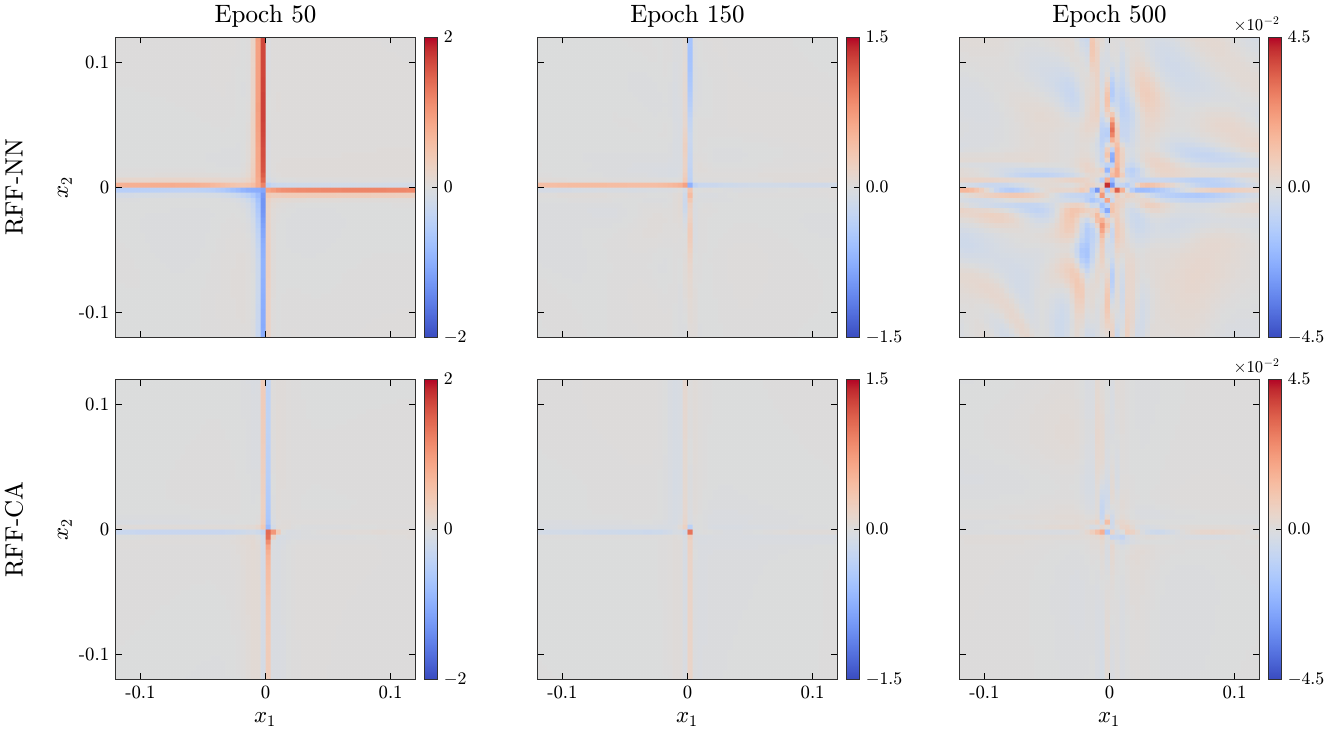}
  \caption{Pointwise error $u(x;\theta)-f_1(x)$ at training-grid points in $[-0.12,0.12]^2$. The upper and lower rows show RFF-NN and RFF-CA, respectively, at epochs $50$, $150$, and $500$ from left to right.}
  \label{fig:regression_visual_checker}
\end{figure}

\begin{figure}[!htbp]
  \centering
  \includegraphics[width=0.95\linewidth]{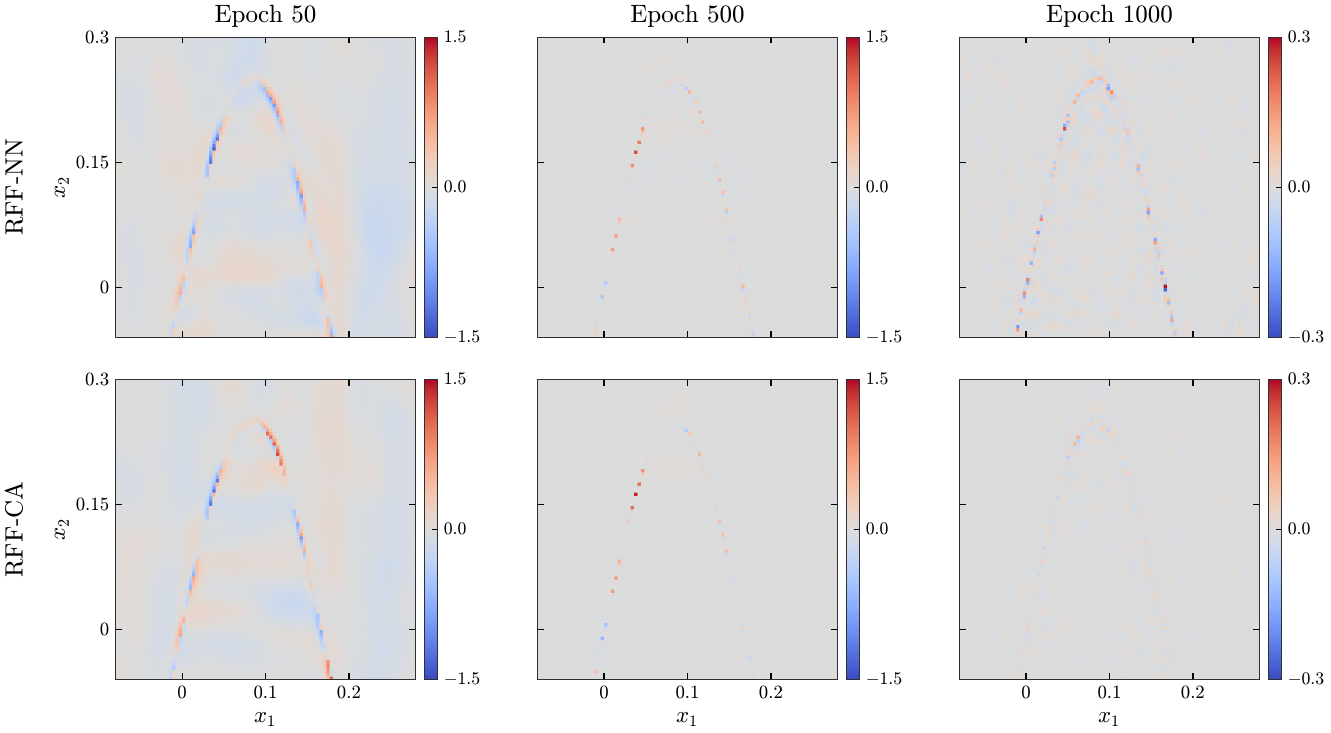}
  \caption{Pointwise error $u(x;\theta)-f_2(x)$ at training-grid points in $[-0.08,0.28]\times[-0.06,0.30]$. The upper and lower rows show RFF-NN and RFF-CA, respectively, at epochs $50$, $500$, and $1000$ from left to right.}
  \label{fig:regression_visual}
\end{figure}

\begin{table}[!htbp]
  \centering
  \caption{Final relative $L^2$ errors on the training grid for the two regression examples. Bold values indicate the smallest error for each function.}
  \label{tab:synthetic_regression_errors}
  \begingroup
  \small
  \setlength{\tabcolsep}{18pt}
  \renewcommand{\arraystretch}{1.15}
  \begin{tabular}{lccc}
    \toprule
    Function & RFF-NN & RFF-CA & NN-CA \\
    \midrule
    $f_1$ & $2.796\times10^{-3}$ & $\boldsymbol{4.394\times10^{-4}}$ & $5.529\times10^{-4}$ \\
    $f_2$ & $1.039\times10^{-2}$ & $\boldsymbol{4.490\times10^{-3}}$ & $5.625\times10^{-3}$ \\
    \bottomrule
  \end{tabular}
  \endgroup
\end{table}
\begin{samepage}
The local error plots in Figs.~\ref{fig:regression_visual_checker} and~\ref{fig:regression_visual} further show smaller error magnitudes for RFF-CA than for RFF-NN in the displayed neighborhoods of both discontinuity interfaces at the final epochs.
\par\end{samepage}

\FloatBarrier
\subsubsection{Image approximation}\label{sec:image_approx}
\begin{figure}[!htbp]
  \centering
  \includegraphics[height=2.8cm]{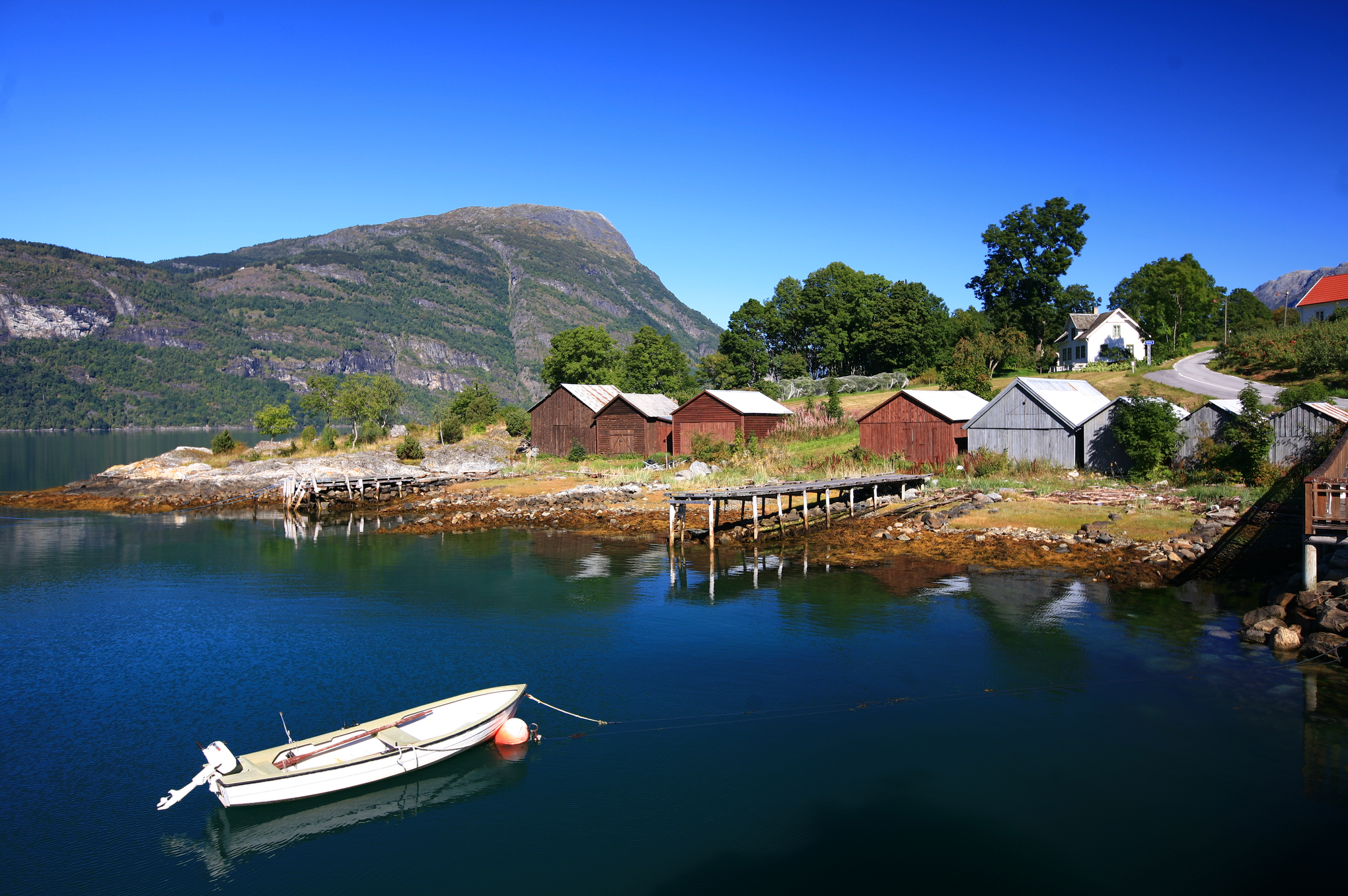}
  \includegraphics[height=2.8cm]{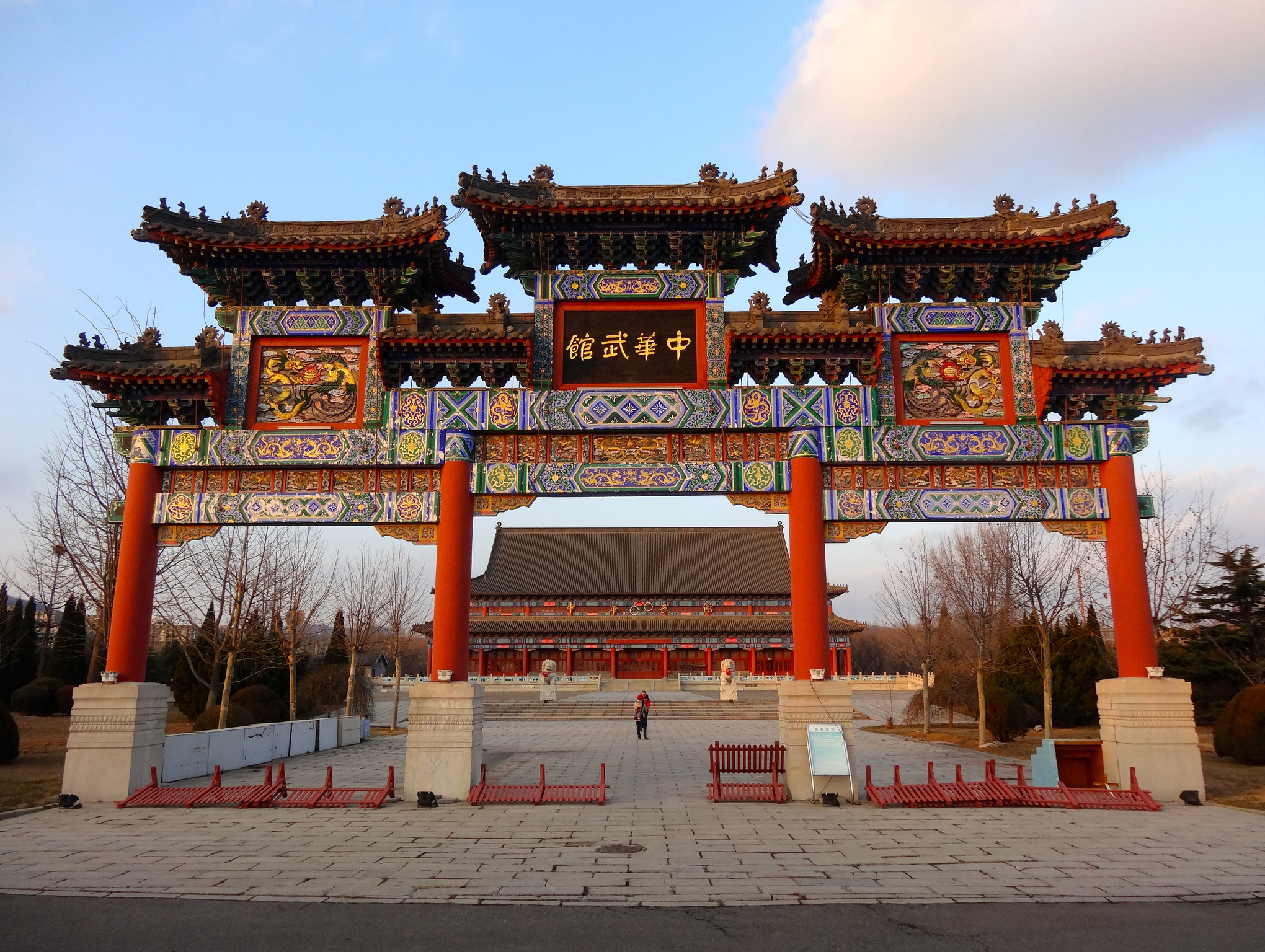}
  \includegraphics[height=2.8cm]{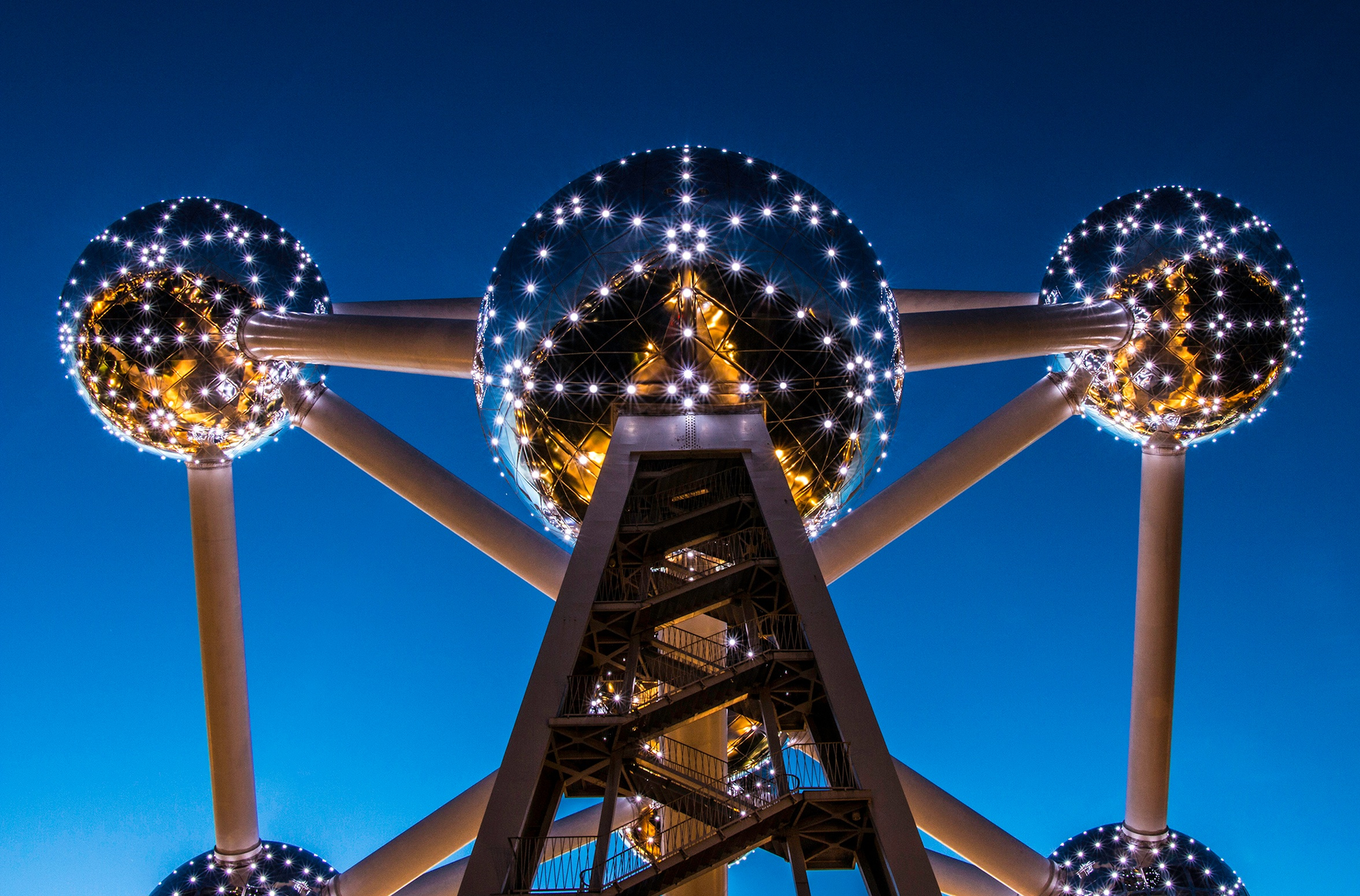}
  \includegraphics[height=2.8cm]{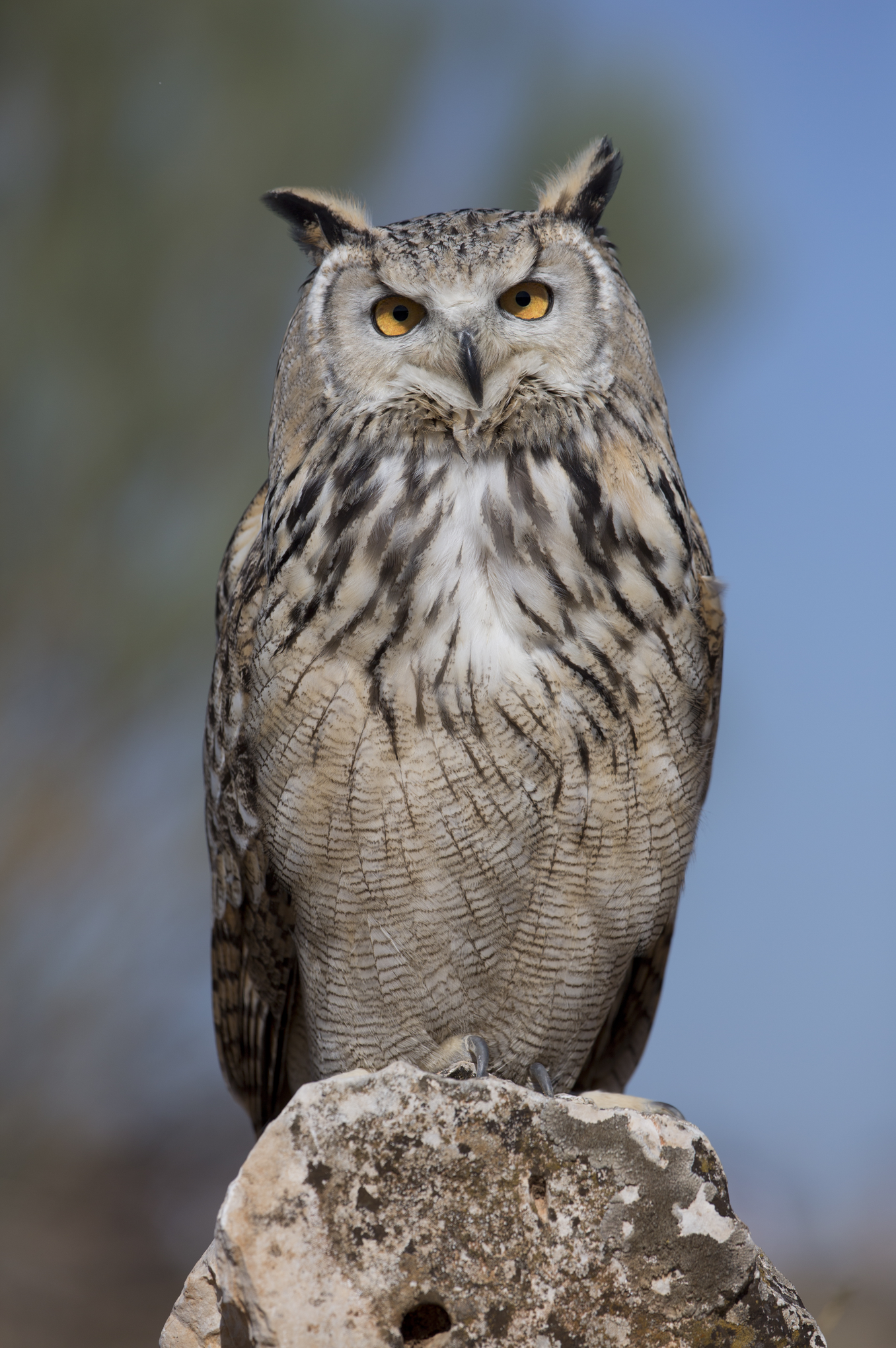}
  \caption{Sample images from the DIV2K validation set used in the regression experiments.
    From left to right, image sizes are
    $2040\times1356$, $2040\times1536$, $2040\times1344$, and $1356\times2040$.}
  \label{fig:div2k-samples}
\end{figure}

We next consider fine-scale textures and sharp edges in the four images from the DIV2K validation set~\cite{agustsson2017ntire} shown in Fig.~\ref{fig:div2k-samples}. Each RGB image is viewed as an $\mathbb{R}^3$-valued function of spatial coordinates, whose components represent the color intensities.

\begingroup\clubpenalty=10000
After downsampling each image by a factor of $4$, we map the pixel centers of the resulting RGB image $I\in[0,1]^{H\times W\times 3}$ to the continuous domain $\Omega=[-1,1]^2$ via
\begin{equation*}
  x_{j} \;=\; \frac{j+\tfrac12}{W}\cdot 2 - 1,
  \qquad
  y_{i} \;=\; \frac{i+\tfrac12}{H}\cdot 2 - 1,
\end{equation*}
for $i=0,\ldots,H-1$ and $j=0,\ldots,W-1$, and define $\xi_{ij}=(x_j,y_i)\in\Omega$.
The regression data consist of the pairs $(\xi_{ij},I_{ij})$, where $I_{ij}\in\mathbb{R}^3$ contains the RGB intensities at pixel $(i,j)$.
We compare RFF-NN, RFF-CA, and NN-CA under this coordinate-to-RGB setting, with a common initial multiscale Fourier feature bank for the three models.
\par\endgroup

For image reconstruction, we report the relative $L^2$ error, the peak signal-to-noise ratio (PSNR), and the relative high-frequency error norm (HFEN).
For the reference image $I$ and its reconstruction $\widehat{I}$, we apply a Laplacian-of-Gaussian (LoG) filter channelwise to obtain
\begin{equation*}
  I^{\mathrm{HP}}=\mathrm{LoG}_\sigma * I,\qquad
  \widehat{I}^{\mathrm{HP}}=\mathrm{LoG}_\sigma * \widehat{I},
\end{equation*}
and define the relative HFEN by
\begin{equation*}
  \mathrm{HFEN}_{\mathrm{rel}}(\widehat{I},I)
  =\frac{\|\,\widehat{I}^{\mathrm{HP}}-I^{\mathrm{HP}}\,\|_2}
  {\|\,I^{\mathrm{HP}}\,\|_2}.
\end{equation*}
Here $*$ denotes discrete convolution with a $15\times15$ LoG kernel and $\sigma\approx1.5$ pixels.

All three models use full-batch Adam with no weight decay and the same learning-rate schedule for all images, decreasing exponentially from $5\times10^{-3}$ to $10^{-4}$ over $5000$ epochs. The corresponding training histories are given in \ref{app:image-regression-curves}.

\begin{table}[H]
  \centering
  \caption{Reconstruction quality after 5000 epochs. Images are numbered from left to right in Fig.~\ref{fig:div2k-samples}. Bold values indicate the best result for each image and metric.}
  \label{tab:img_results}
  \begingroup
  \small
  \setlength{\tabcolsep}{3pt}
  \renewcommand{\arraystretch}{1.2}
  \begin{tabular*}{\linewidth}{@{\extracolsep{\fill}}l*{9}{r}@{}}
    \toprule
    & \multicolumn{3}{c}{Relative $L^2$ error $\downarrow$}
    & \multicolumn{3}{c}{PSNR (dB) $\uparrow$}
    & \multicolumn{3}{c}{Relative HFEN $\downarrow$} \\
    \cmidrule(lr){2-4}\cmidrule(lr){5-7}\cmidrule(lr){8-10}
    Image & RFF-NN & RFF-CA & NN-CA & RFF-NN & RFF-CA & NN-CA & RFF-NN & RFF-CA & NN-CA \\
    \midrule
    1 & 0.0878 & \textbf{0.0632} & 0.0640 & 27.83 & \textbf{30.68} & 30.58 & 0.3746 & \textbf{0.2278} & 0.2363 \\
    2 & 0.0992 & \textbf{0.0775} & 0.0776 & 24.92 & \textbf{27.07} & 27.04 & 0.4656 & \textbf{0.3195} & 0.3202 \\
    3 & 0.1108 & 0.0721 & \textbf{0.0658} & 27.32 & 31.05 & \textbf{31.84} & 0.3386 & 0.1774 & \textbf{0.1570} \\
    4 & 0.0680 & \textbf{0.0536} & 0.0548 & 29.75 & \textbf{31.82} & 31.61 & 0.3244 & \textbf{0.2251} & 0.2319 \\
    \bottomrule
  \end{tabular*}
  \endgroup
\end{table}

\begin{figure}[H]
  \centering
  \includegraphics[width=0.9\linewidth]{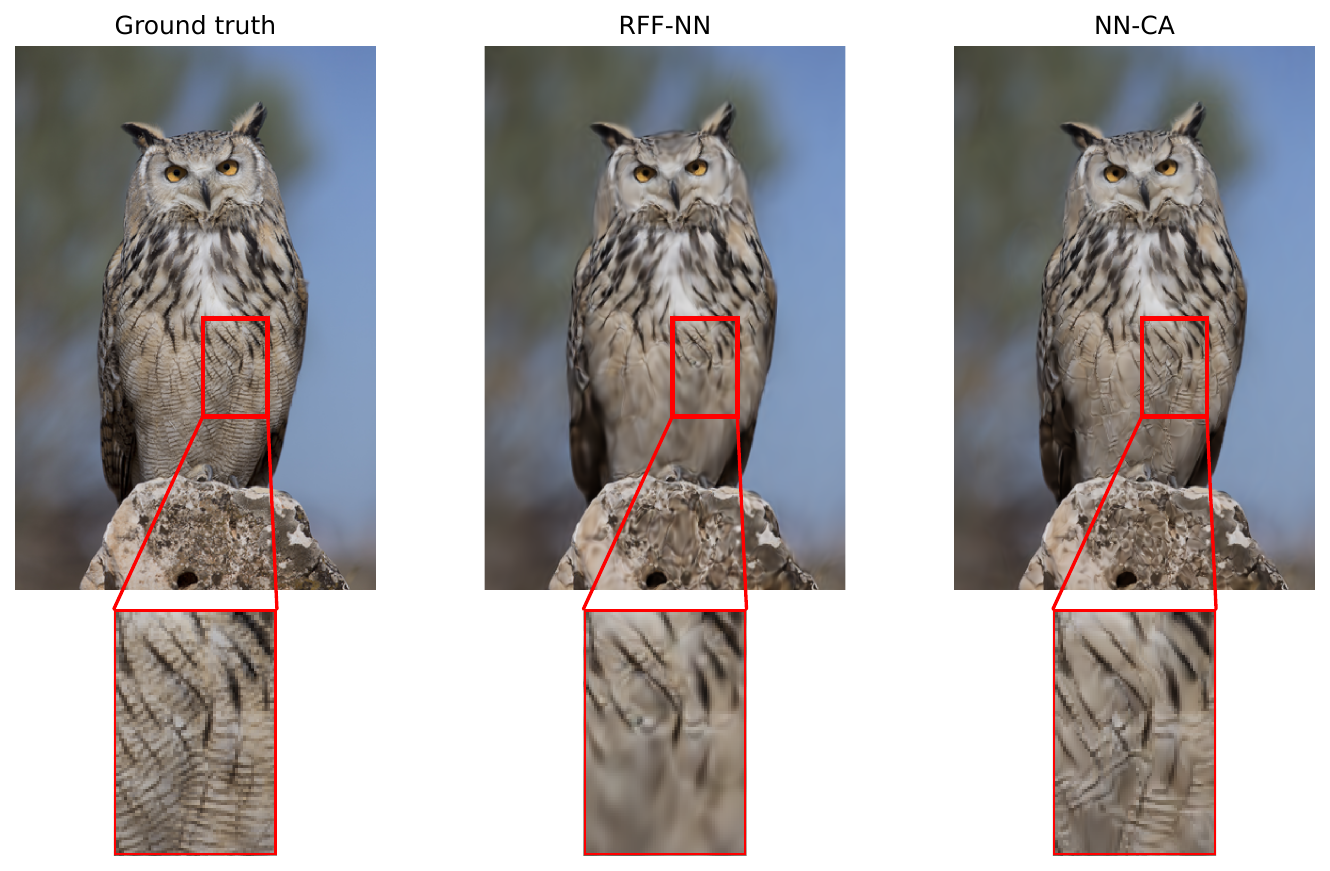}
  \caption{The reference image (left) and the RFF-NN (middle) and NN-CA (right) reconstructions after $5000$ epochs are shown with enlarged views of the boxed regions below.}
  \label{fig:img_visual}
\end{figure}
For all four images, both RFF-CA and NN-CA achieve lower relative $L^2$ error and HFEN and higher PSNR than RFF-NN after $5000$ epochs (Table~\ref{tab:img_results}). RFF-CA gives the best results for images 1, 2, and 4, whereas NN-CA performs best for image 3. In the zoomed regions of the owl image (Fig.~\ref{fig:img_visual}), NN-CA preserves fine-scale textures and edge contrast more faithfully than RFF-NN, whose reconstructions appear smoother and contain fewer high-frequency details.

\FloatBarrier
\subsubsection{Adaptive frequency enhancement}
We study the effect of AFE on approximation accuracy using the one-dimensional periodic function on $\Omega=[0,1)$ defined by
\begin{equation}
  u(x)
  =
  \sin(2\pi\cdot 2x+0.1)
  +0.5\,\sin(2\pi\cdot 20x+0.3)
  +0.25\,\cos(2\pi\cdot 40x-0.2).
  \label{eq:afe-1d-target}
\end{equation}
Its sparse spectrum allows the extracted frequencies to be compared with the known Fourier modes.

We use RFF-CA with a single-scale bank ($K=0$), $\sigma_\omega=1$, and $L=3$ residual blocks. Training and error evaluation use uniform periodic grids of $8188$ points each, with the evaluation grid shifted by half a mesh spacing relative to the training grid. Optimization uses full-batch Adam with an initial learning rate of $10^{-3}$, multiplied by $0.85$ every $500$ epochs.

\begin{figure}[!htbp]
  \centering
  \includegraphics[width=0.48\linewidth]{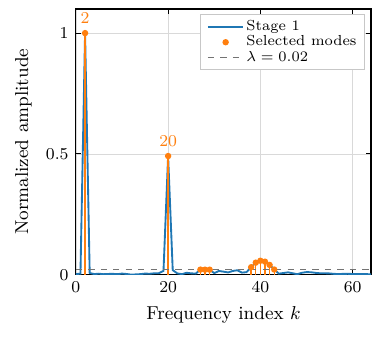}

  \caption{Normalized DFT spectrum of the Stage~1 prediction, with orange markers indicating the selected modes and the dashed line marking the threshold $\lambda=0.02$.}
  \label{fig:afe-frequency}
\end{figure}

After $5000$ epochs (Stage~1), we apply the adaptive frequency selection procedure in Section~\ref{subsec:afe} to the intermediate approximation, using a uniform DFT grid of $4096$ points and a threshold of $\lambda=0.02$. The selected frequencies include the three modes $k=2,20,40$ of the reference function, as shown in Fig.~\ref{fig:afe-frequency}. Starting from the same network parameters and Adam optimizer state, we then continue training for another $5000$ epochs (Stage~2) using either the original feature bank or the bank enriched by AFE, with the learning-rate schedule continued in both cases. For AFE, we set the attention-mask parameter to $\eta=-6$.

\begin{figure}[!htbp]
  \centering
  \includegraphics[width=\linewidth]{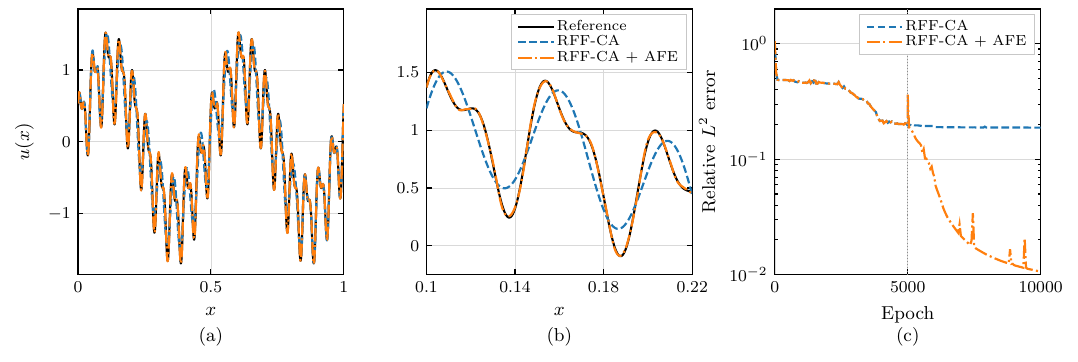}

  \caption{Effect of AFE on the RFF-CA approximation, shown through the final approximations (a), a local view on $[0.10,0.22]$ (b), and the relative $L^2$ error histories (c). The vertical dotted line marks the start of Stage~2, and the AFE curve includes the error immediately after the feature bank is augmented.}
  \label{fig:afe-results}
\end{figure}

After the same total of $10{,}000$ training epochs, the final relative $L^2$ error on the evaluation grid decreases from $1.88\times10^{-1}$ for the original bank to $1.06\times10^{-2}$ for AFE, a reduction by a factor of approximately $17.7$ in this example (Fig.~\ref{fig:afe-results}). Although the error initially increases when the feature bank is augmented, it decreases overall during Stage~2.

\FloatBarrier
\subsection{PDE problems}
We next assess the two-network formulation introduced in Section~\ref{sec:PDE}, beginning with a one-dimensional Poisson example to compare mixing strategies, network architectures, and amplitude weighting. The subsequent examples consider spatially varying oscillations in two dimensions and a discontinuous reaction coefficient in three dimensions.
\subsubsection{1D Poisson equation}
Consider the Poisson equation in $\Omega=[-1,1]$ with Dirichlet boundary conditions,
\begin{align*}
  -\Delta u(x) = f(x),
\end{align*}
for which the exact solution is chosen as
\begin{equation*}
  u(x)=\sin(0.1\pi x)+0.2\sin(\pi x)+0.4\sin\big((\nu/3)\pi x\big)
  +0.6\sin\big(2(\nu/3)\pi x\big)+\sin(\nu\pi x),
\end{equation*}
where $\nu=100$, so that $u$ contains a mixture of very low- and very
high-frequency components. The source term and boundary values are determined by this exact solution.

For the one-dimensional experiments, the low-frequency network has two hidden layers of width $64$ with $\tanh$ activation. Unless specified below, we set $\sigma_\omega=0.1$, $\beta_0=\log(1+e^{-5})$, and $\gamma=10^3$. We train for $20{,}000$ epochs using AdamW with weight decay $10^{-2}$ and no gradient clipping, starting from a learning rate of $10^{-3}$ that is halved every $2000$ epochs. Each epoch uses $10^4$ uniformly sampled interior points and $20$ boundary samples, ten at each endpoint, while relative $L^2$ errors are evaluated on a uniform grid of $2001$ points.

With NN-CA as the high-frequency component $u_h$, we first compare four choices of the mixing coefficient $\alpha$, namely fixed values of $0$ and $1$, a learnable coefficient initialized at $1$, and a coefficient computed analytically at each epoch. The fixed and learnable strategies use the loss $\mathcal{L}$, whereas the analytical strategy computes $\alpha$ from $\mathcal{L}$ and trains the network parameters using $\widehat{\mathcal{L}}$, as described in Section~\ref{sec:PDE}.

\begin{figure}[!htbp]
  \centering
  \includegraphics[width=0.88\linewidth]{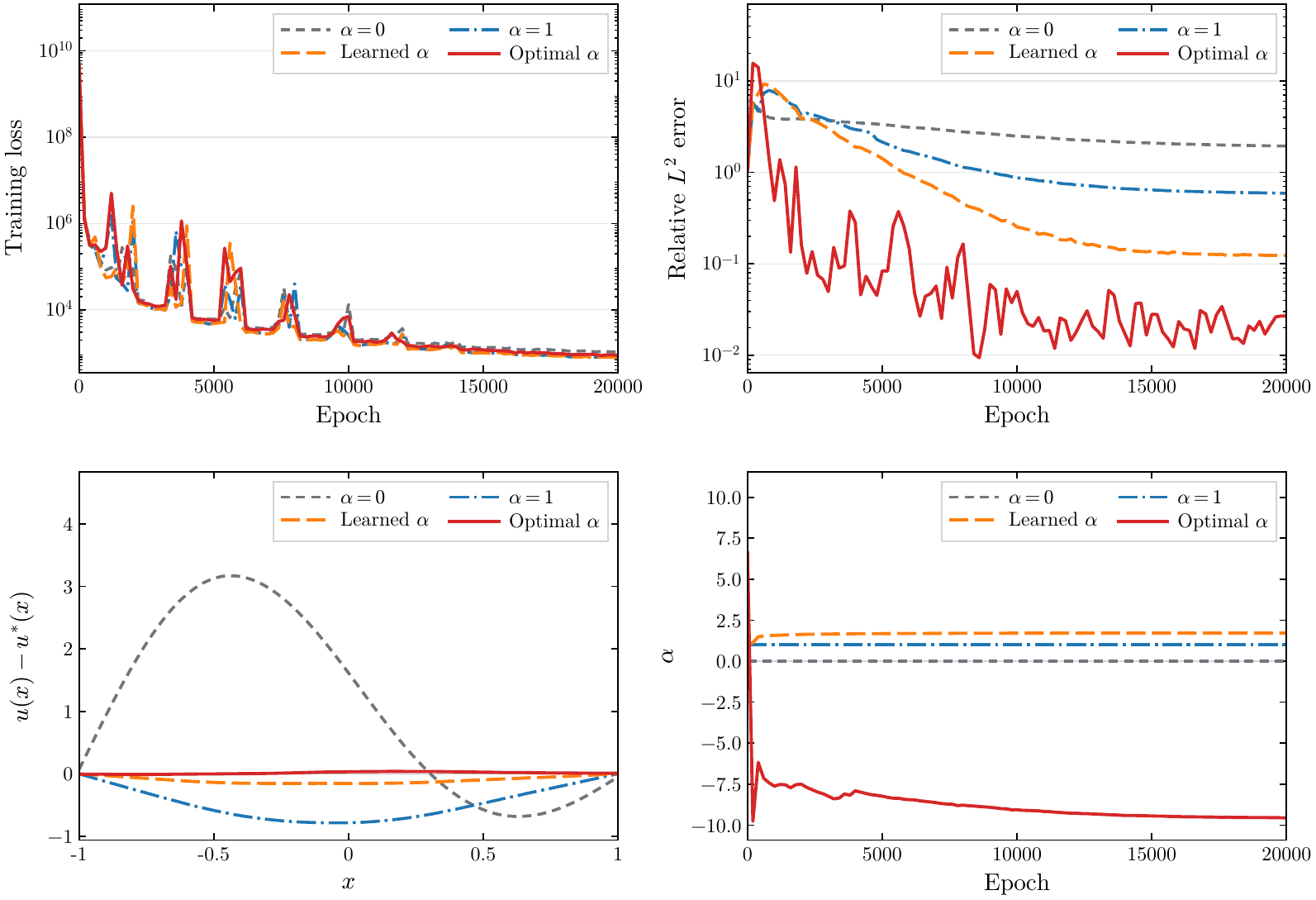}

  \caption{For the four mixing strategies in the 1D Poisson example, the upper panels compare the training losses (left) and relative $L^2$ error histories (right). The lower panels show the final pointwise error, computed as the prediction minus the exact solution (left), and the evolution of $\alpha$ (right).}
  \label{fig:poisson_benchmark_error_loss_evolution}
\end{figure}

\begingroup\clubpenalty=10000
In Fig.~\ref{fig:poisson_benchmark_error_loss_evolution}, the analytical strategy achieves the smallest final relative $L^2$ error, $2.69\times10^{-2}$, compared with $1.24\times10^{-1}$ for the learnable coefficient. The training-loss curves describe the evolution of the respective objectives, while the relative $L^2$ errors provide a common measure of approximation accuracy.
Under the reported training configuration, the single-network case
  $\alpha=0$ yields a larger approximation error than the two-network alternatives.
  During the initial training phase, the analytically computed coefficient
  takes large-magnitude negative values, and the relative $L^2$ error
  decreases rapidly.
\par\endgroup

\begin{figure}[!htbp]
  \centering
  \includegraphics[width=0.95\linewidth]{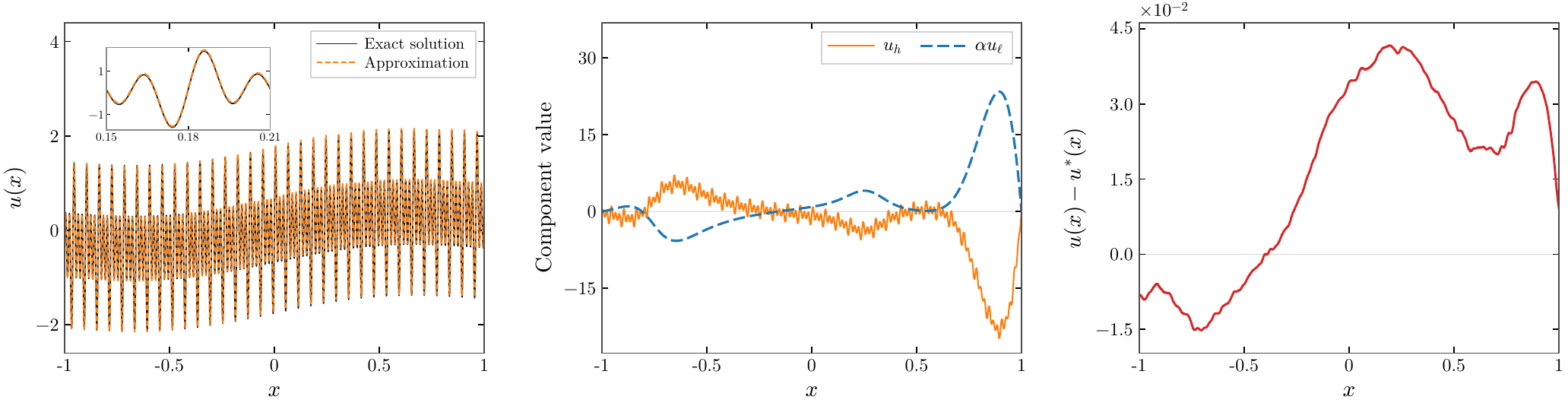}

  \caption{The prediction with the analytically computed mixing coefficient is compared with the exact solution in the left panel, including a local view on $[0.15,0.21]$. The middle panel shows the components $u_h$ and $\alpha u_l$, and the right panel gives the pointwise error.}
  \label{fig:poisson_benchmark_optimal}
  \par\vspace{10pt}
  \includegraphics[width=0.9\linewidth]{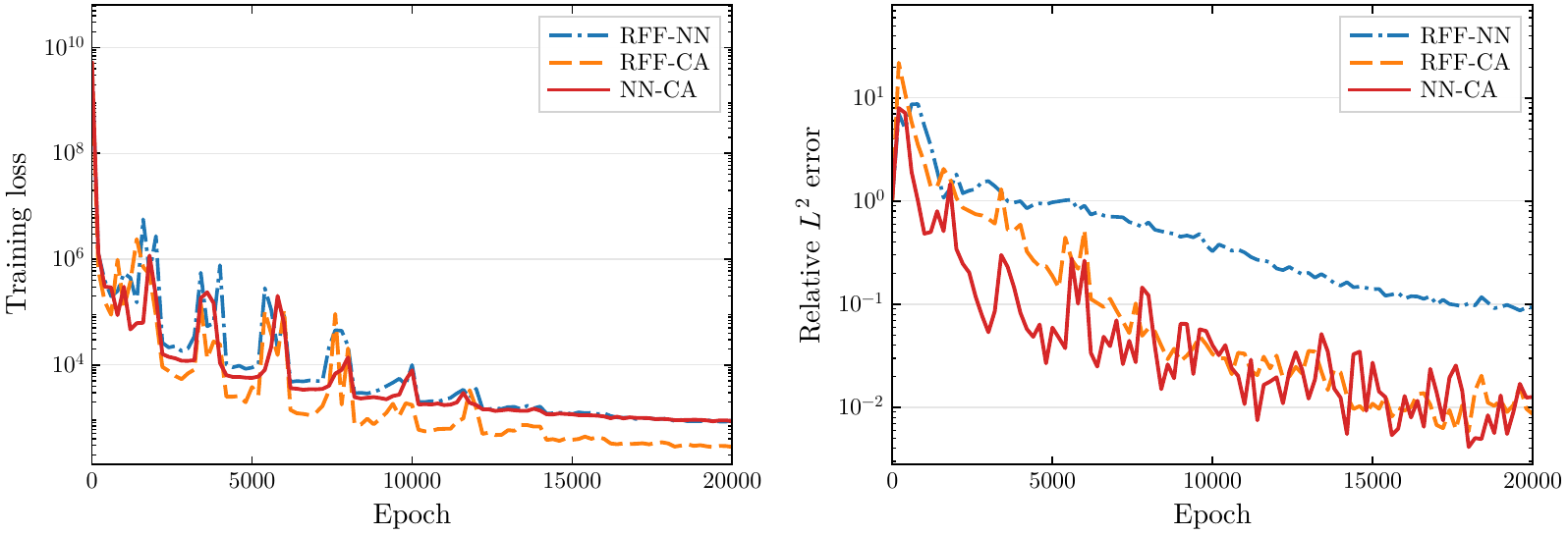}

  \caption{Training histories of RFF-NN, RFF-CA, and NN-CA for the 1D Poisson problem with the analytically computed mixing coefficient, showing the loss $\widehat{\mathcal{L}}$ (left) and relative $L^2$ error (right).}
  \label{fig:poisson_benchmark_structure}
  \par\vspace{10pt}
  \includegraphics[width=0.9\linewidth]{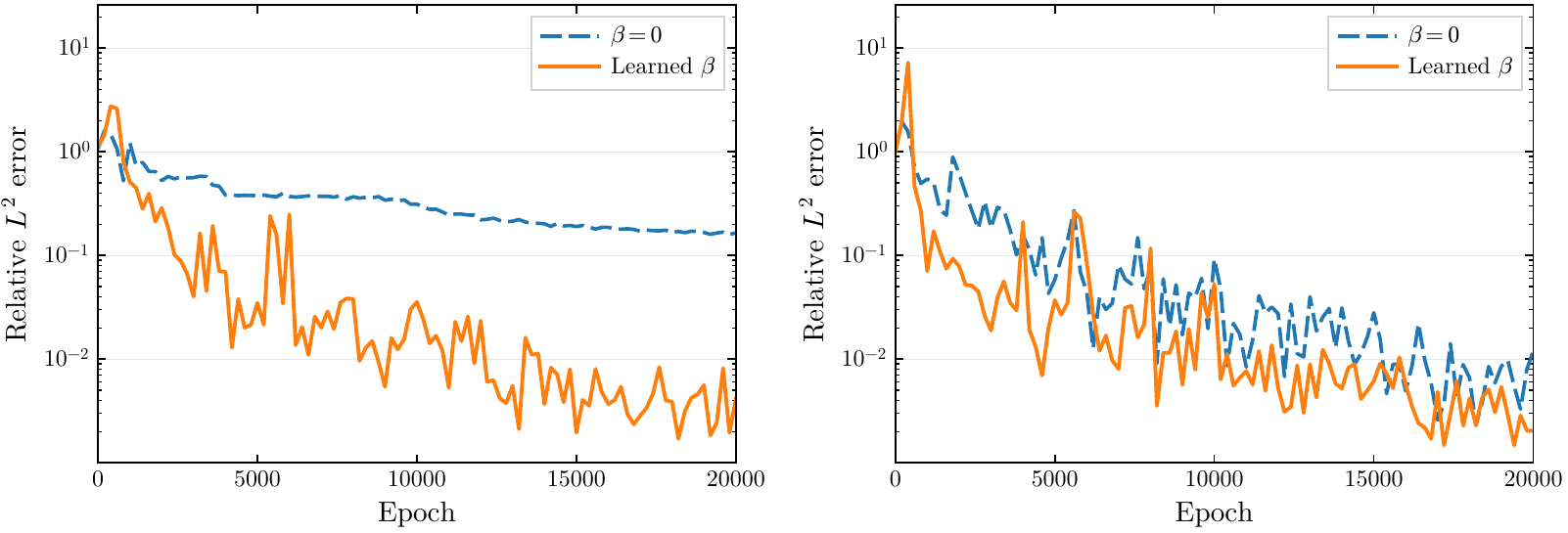}

  \caption{The relative $L^2$ error histories compare RFF-CA with a learnable envelope parameter $\beta$ in its multiscale RFF bank against the equal-amplitude variant ($\beta\equiv0$), for boundary-penalty weights $\gamma=10^3$ (left) and $\gamma=10^4$ (right).}
  \label{fig:poisson_benchmark_scale}
\end{figure}

The component plots in Fig.~\ref{fig:poisson_benchmark_optimal} show that $u_h$ captures the oscillatory structure, while $\alpha u_l$ provides a smoother correction, together approximating the reference solution with small pointwise error.

Next, we compare RFF-NN, RFF-CA, and NN-CA as choices for $u_h$, using the analytically computed mixing coefficient and the same initial multiscale Fourier feature bank. As shown in Fig.~\ref{fig:poisson_benchmark_structure}, RFF-CA and NN-CA attain smaller final relative $L^2$ errors than RFF-NN, with RFF-CA giving the most accurate approximation in this comparison.

Finally, for a broader frequency distribution with $\sigma_\omega=0.02$, we compare RFF-CA with a learnable amplitude envelope against its equal-amplitude variant ($\beta\equiv0$). Both use the analytically computed mixing coefficient and the loss $\widehat{\mathcal{L}}$, with boundary-penalty weights $\gamma=10^3$ and $\gamma=10^4$. As shown in Fig.~\ref{fig:poisson_benchmark_scale}, learning the amplitude envelope gives a smaller final relative $L^2$ error for both values of $\gamma$, while increasing $\gamma$ also improves the final accuracy of the equal-amplitude variant.

\FloatBarrier
\subsubsection{2D Poisson equation}
To assess the two-network approximation of a solution with spatially varying oscillations in two dimensions, we consider the Poisson equation in $\Omega=[-1,1]^2$,
\begin{align*}
  -\Delta u(x) = f(x), \qquad x=(x_1,x_2)\in\Omega,
\end{align*}
with Dirichlet boundary conditions.
We choose the exact solution
\begin{equation*}
  u(x)=\sin(\mu x_1^2)+\sin(\mu x_2^2),
\end{equation*}
We consider $\mu=50$ and $100$, with the latter corresponding to a more oscillatory solution.
The source term and boundary values are determined by this exact solution.

We compare NN, RFF-NN, and RFF-CA as choices for the high-frequency component $u_h$. The low-frequency component $u_l$ is represented by a network with two hidden layers of width $64$ and $\tanh$ activation. The NN variant takes the spatial coordinates as input and retains the residual blocks and output layers of RFF-NN. All three models use the analytically computed mixing coefficient and the modified residual-based loss $\widehat{\mathcal{L}}$, with boundary-penalty weight $\gamma=10^4$. For the two RFF-based models, we set $\sigma_\omega=0.1$ and $\beta_0=\log(1+e^{-5})$.

At each epoch, we uniformly sample $10^4$ interior collocation points and $4000$ boundary points. The models are trained using SOAP~\cite{vyas2024soap,khodakarami2026spectral} for $10{,}000$ epochs, starting from a learning rate of $10^{-3}$ that is halved every $2000$ epochs. For each value of $\mu$, the three models use the same initial low-frequency network parameters and collocation samples, and the two RFF-based models additionally share the initial Fourier feature bank.

\begin{figure}[!htbp]
  \centering
  \includegraphics[width=0.95\linewidth]{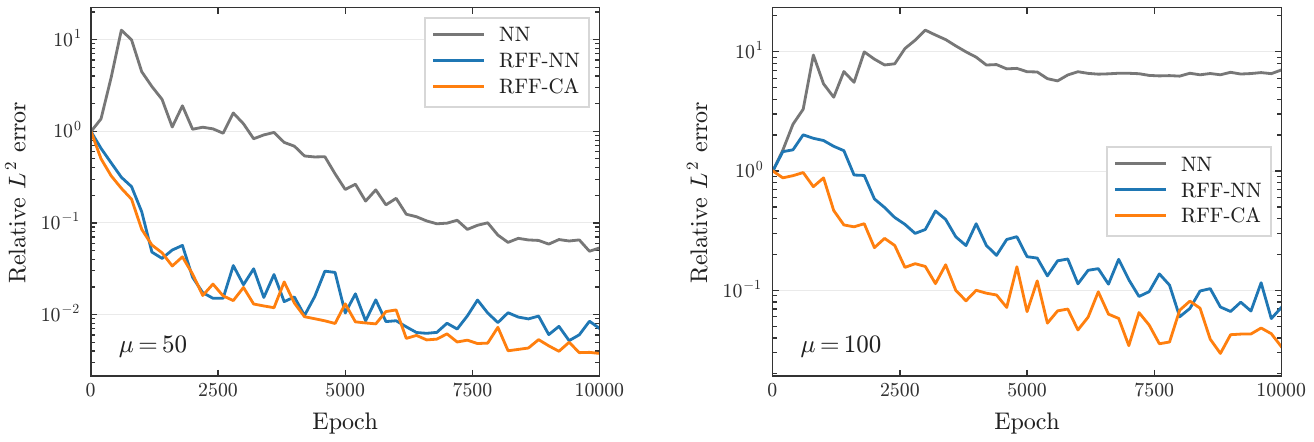}

  \caption{The relative $L^2$ errors of NN, RFF-NN, and RFF-CA for $\mu=50$ (left) and $\mu=100$ (right) are plotted over the first $10{,}000$ epochs of SOAP training.}
  \label{fig:broad_frequency_mu}
\end{figure}

Figure~\ref{fig:broad_frequency_mu} shows that RFF-CA attains the smallest relative $L^2$ error for both values of $\mu$. For $\mu=50$, NN, RFF-NN, and RFF-CA achieve errors of $5.34\times10^{-2}$, $7.03\times10^{-3}$, and $3.77\times10^{-3}$, respectively. When $\mu$ increases to $100$, the error of NN rises to $7.03$, whereas RFF-NN and RFF-CA attain errors of $7.21\times10^{-2}$ and $3.37\times10^{-2}$. Taken together, these experiments suggest that both the representation and the optimization strategy are important for accurately approximating highly oscillatory solutions.

\begin{figure}[!htbp]
  \centering
  \includegraphics[width=\linewidth]{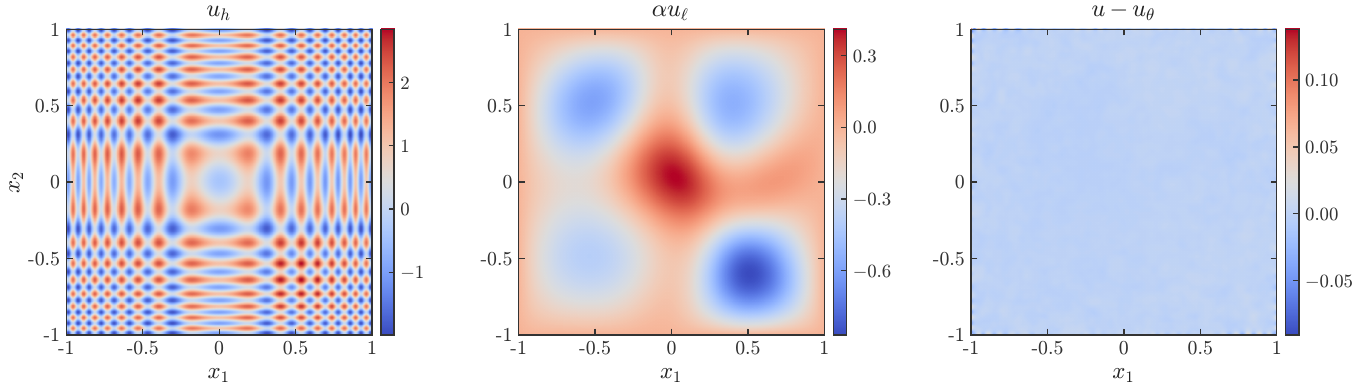}

  \caption{The left and middle panels show the RFF-CA components $u_h$ and $\alpha u_l$, respectively, for $\mu=50$ after $10{,}000$ epochs of SOAP training, while the right panel shows the pointwise error defined between the exact solution and the prediction.}
  \label{fig:broad_frequency_components}
\end{figure}

For the RFF-CA approximation with $\mu=50$ after $10{,}000$ epochs, Fig.~\ref{fig:broad_frequency_components} shows that $u_h$ represents the oscillatory structure, whereas $\alpha u_l$ provides a smoother contribution.

\FloatBarrier
\subsubsection{A three-dimensional Poisson--Boltzmann equation}
\begin{figure}[!htbp]
  \centering
  \includegraphics[height=5cm]{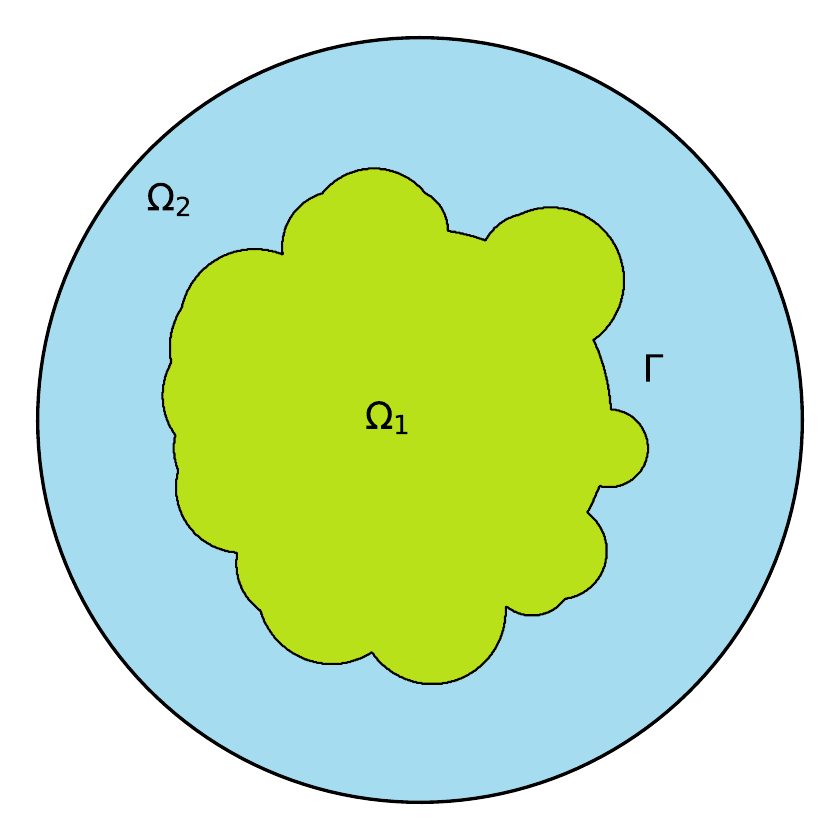}
  \includegraphics[height=5cm]{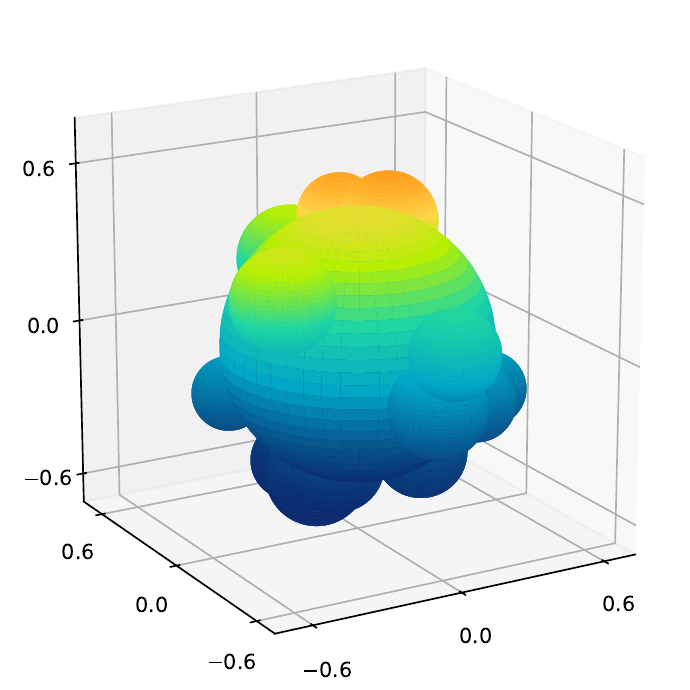}
  \caption{The left panel shows a schematic of the coefficient partition into two subdomains and their internal interface, while the right panel illustrates the three-dimensional subdomain $\Omega_1$ formed by the union of overlapping balls.}
  \label{fig:singular_domain}
\end{figure}

We finally consider a linear Poisson--Boltzmann equation \cite{liu2020multi} on the unit ball $\Omega=\{x\in\mathbb{R}^3:|x|<1\}$,
\begin{align*}
  -\nabla\cdot (\epsilon (x)\nabla u(x)) + \kappa(x)u(x) = f(x), \quad x\in \Omega,
\end{align*}
with homogeneous Dirichlet boundary conditions
\begin{align*}
  u|_{\partial \Omega}=0.
\end{align*}
The diffusion coefficient is $\epsilon(x)\equiv1$, and the reaction coefficient is
\begin{equation*}
  \kappa(x)=\begin{cases}
    1, & x\in\Omega_1, \\
    5, & x\in\Omega_2,
  \end{cases}
  \qquad \Omega_2=\Omega\setminus\Omega_1.
\end{equation*}
Here $\Omega_1$ is the union of a ball of radius $0.5$ centered at the origin and $20$ smaller balls. Their centers are sampled uniformly on the surface of the central ball, and their radii are sampled uniformly from $[0.1,0.2]$. The intersections of the balls give the internal interface $\Gamma=\partial\Omega_1$ a nonsmooth geometry, as illustrated in Fig.~\ref{fig:singular_domain}.

We choose $f=-\Delta u+\kappa u$ so that the exact solution is
\begin{align*}
  u(x) = \frac{e^{\sin \mu x_1 + \sin \mu x_2 + \sin \mu x_3}}{|x|^2 + 1}\left(|x|^2-1\right),\qquad \mu=15.
\end{align*}
This smooth, oscillatory solution satisfies the boundary condition and the transmission conditions $[u]=0$ and $[\epsilon\,\partial u/\partial n]=0$ across $\Gamma$, where $[\cdot]$ denotes the jump across the interface.

For this example, we train the two-network approximation using the Deep Ritz variational formulation, with loss
\begin{align*}
  \mathcal{L}_{\mathrm{Ritz}}(u(\cdot;\theta))
  &= \frac{1}{2}\int_\Omega \Bigl(\epsilon(x)|\nabla u(x;\theta)|^2 + \kappa(x)\,u(x;\theta)^2\Bigr)\,\mathrm{d}x \\
  &\quad -\int_\Omega f(x)\,u(x;\theta)\,\mathrm{d}x
  +\gamma\int_{\partial\Omega} |u(x;\theta)-g(x)|^2\,\mathrm{d}S,
\end{align*}
The boundary penalty uses $g\equiv0$ and $\gamma=10^4$.

We compare RFF-NN and RFF-CA as choices for $u_h$, with $\sigma_\omega=1$ and $\beta_0=\log(1+e)$. The low-frequency network $u_l$ has three hidden layers of width $128$ with $\tanh$ activation, and the mixing coefficient $\alpha$ is initialized at $1$ and optimized jointly with the network weights using the same Ritz loss.

The models are trained for $10{,}000$ epochs using AdamW with weight decay $10^{-2}$ and gradient clipping at $10$, starting from a learning rate of $10^{-3}$ that is multiplied by $0.6$ every $1000$ epochs. At each epoch, we uniformly sample $5000$ interior points and $4000$ boundary points, while the relative $L^2$ error is evaluated on a fixed set of $20{,}000$ points sampled uniformly in the three-dimensional domain.

\begin{figure}[!htbp]
  \centering
  \includegraphics[width=0.95\linewidth]{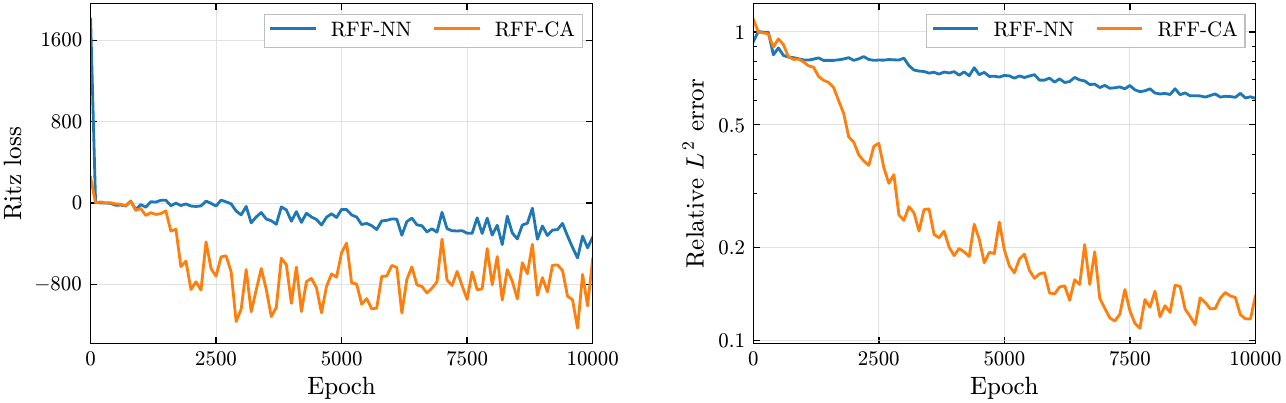}
  \caption{Training Ritz loss (left) and relative $L^2$ error on a fixed set of three-dimensional evaluation points (right) for RFF-NN and RFF-CA.}
  \label{fig:singular_results}
\end{figure}

\begin{figure}[!htbp]
  \centering
  \includegraphics[width=\linewidth]{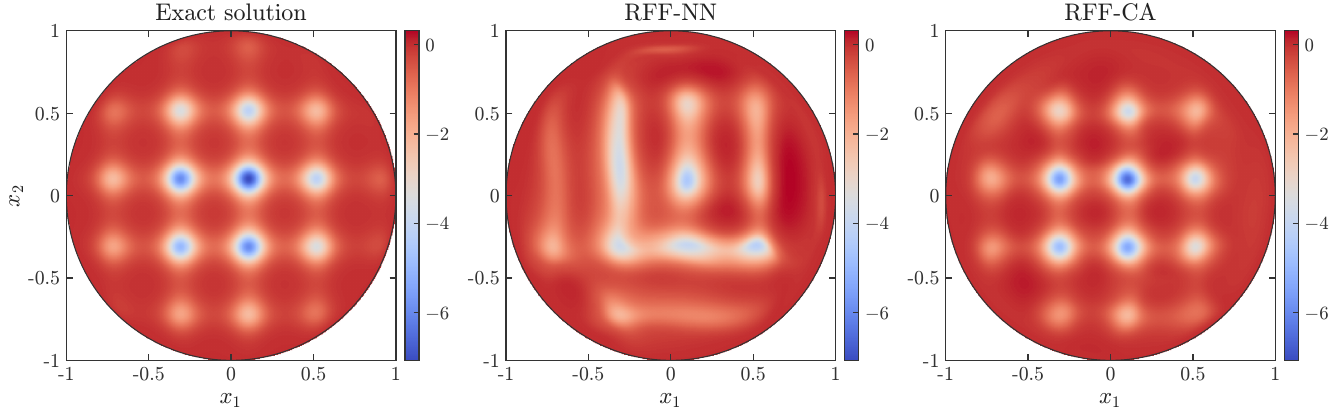}
  \caption{The exact solution (left) and the RFF-NN (middle) and RFF-CA (right) approximations after $10{,}000$ training epochs are shown on the plane $x_3=0$.}
  \label{fig:singular_prediction}
\end{figure}

The training histories in Fig.~\ref{fig:singular_results} show that RFF-CA attains a relative $L^2$ error of $0.139$ after $10{,}000$ epochs, compared with $0.610$ for RFF-NN. In the solution slices on the plane $x_3=0$ shown in Fig.~\ref{fig:singular_prediction}, RFF-CA also approximates the oscillatory structure more accurately.

\FloatBarrier
\section{Conclusion}
\label{sec:conclusion}
In this work, we introduced a cross-attention-based framework to alleviate spectral bias in
high-frequency function approximation and PDE learning. By combining a multiscale random Fourier feature bank with learnable amplitude weighting and cross-attention, the proposed method provides an input-dependent mechanism for emphasizing informative scales. We further developed an adaptive frequency enhancement strategy that enriches the Fourier representation with dominant modes identified through discrete Fourier analysis of intermediate approximations. For PDE problems, we introduced a two-network formulation to improve the approximation of low- and high-frequency solution components, with a learnable or analytically computed mixing coefficient balancing the contributions of the two subnetworks.
However, several important issues remain. A rigorous analysis of approximation and
optimization dynamics for cross-attention-based multiscale Fourier representations is
still lacking. In addition, the current DFT-guided enhancement is most natural for periodic
or grid-friendly settings; extending it to nonperiodic and geometrically complex domains and to
  higher dimensions deserves further investigation. Finally, integrating the proposed spectral
control mechanism into broader physics-informed operator learning methods is a promising
direction for future research.

\section*{Acknowledgements}
Xiaodong Feng is supported by the Guangdong and Hong Kong Universities ``1+1+1'' Joint
Research Collaboration Scheme (Project No.~2025A0505000014) and the Start-up Research
Fund of Beijing Normal-Hong Kong Baptist University.
Tao Zhou is supported by the Science Challenge Project TZ2025006, the NSF of China
under grants 12288201 and 12461160275, and the Strategic Priority Research Program
of the Chinese Academy of Sciences under grant XDA0480504.

\appendix
\section{A mode-wise linearized analysis of differential-loss weighting}
\label{appendix:1}
This appendix relates the frequency weights introduced by differential losses to the decay of approximation errors under parameter gradient flow.

Consider $-\partial_{xx}u=f$ on $\Omega=[-1,1]$, and assume that $u,u(\cdot;\theta)\in H^2(\Omega)\cap H_0^1(\Omega)$. With $e(x;\theta)=u(x;\theta)-u(x)$, define
\begin{align}
  \mathcal{L}_0(\theta)
  &=\frac{1}{2}\|e(\cdot;\theta)\|_{L^2(\Omega)}^2,
  \label{eq:appendix-L0}\\
  \mathcal{L}_1(\theta)
  &=\frac{1}{2}\|\partial_xe(\cdot;\theta)\|_{L^2(\Omega)}^2,
  \label{eq:appendix-L1}\\
  \mathcal{L}_2(\theta)
  &=\frac{1}{2}\|\partial_{xx}e(\cdot;\theta)\|_{L^2(\Omega)}^2.
  \label{eq:appendix-L2}
\end{align}
Here $\mathcal{L}_0$ is the regression loss. For the Deep Ritz energy
\begin{equation*}
  J(v)=\frac{1}{2}\int_\Omega|\partial_xv|^2\,\mathrm{d}x
  -\int_\Omega fv\,\mathrm{d}x,\qquad v\in H_0^1(\Omega),
\end{equation*}
integration by parts and the shared homogeneous boundary conditions give $J(u(\cdot;\theta))-J(u)=\mathcal{L}_1(\theta)$. Since $-\partial_{xx}u(\cdot;\theta)-f=-\partial_{xx}e(\cdot;\theta)$, the functional $\mathcal{L}_2$ equals half the squared $L^2$ norm of the Poisson residual. Boundary penalties are omitted to isolate the interior frequency weighting.

Let $\varphi_j(x)=\sin(j\pi x)$ and consider
\begin{equation*}
  u(x;\theta)=c_1(\theta)\varphi_1(x)+c_k(\theta)\varphi_k(x),
  \qquad
  u(x)=\varphi_1(x)+a\varphi_k(x),
\end{equation*}
where $k\geq2$ is an integer, $a\in\mathbb{R}$, and $c_1,c_k$ are continuously differentiable functions of $\theta$. Set $e_1=c_1-1$, $e_k=c_k-a$, and $\mathbf{e}=(e_1,e_k)^\top$. Orthogonality of the Fourier modes and their derivatives yields
\begin{equation}
  \mathcal{L}_m(\theta)=\frac{1}{2}\mathbf{e}^{\top}W_m\mathbf{e},
  \qquad
  W_m=\operatorname{diag}\bigl(\pi^{2m},(k\pi)^{2m}\bigr),
  \qquad m=0,1,2.
  \label{eq:appendix-modal-loss}
\end{equation}
Under parameter gradient flow $\dot{\theta}=-\nabla_\theta\mathcal{L}_m$, the identity $\dot e_j=\langle\nabla_\theta c_j,\dot\theta\rangle$ gives
\begin{equation}
  \frac{\mathrm{d}\mathbf{e}}{\mathrm{d}t}=-G_\theta W_m\mathbf{e},
  \qquad
  G_\theta=\bigl[\langle\nabla_\theta c_i,\nabla_\theta c_j\rangle\bigr]_{i,j\in\{1,k\}}.
  \label{eq:appendix-general-dynamics}
\end{equation}
The symmetric positive semidefinite matrix $G_\theta$ is the restriction of the tangent kernel~\cite{jacot2018neural} to the two-mode subspace. It describes the parameter sensitivity and coupling of the modal coefficients, whereas $W_m$ contains the weights imposed by the loss.

To compare the objectives at the same linearization, we replace $G_\theta$ in \eqref{eq:appendix-general-dynamics} by $\overline G=G_{\theta(0)}$, with the same initial error. For this linearized model, the weighted error $\mathbf z=W_m^{1/2}\mathbf e$ satisfies
\begin{equation*}
  \dot{\mathbf z}=-S_m\mathbf z,
  \qquad S_m=W_m^{1/2}\overline G W_m^{1/2}.
\end{equation*}
Since $S_m$ is symmetric and positive semidefinite, write $S_m=U_m\Lambda_mU_m^\top$, where $U_m$ is orthogonal and $\Lambda_m=\operatorname{diag}(\mu_{m,1},\mu_{m,2})$ with $\mu_{m,i}\geq0$. Its spectral decomposition gives
\begin{equation*}
  \mathbf z(t)=U_m\exp(-\Lambda_m t)U_m^\top\mathbf z(0),
  \qquad
  \mathbf e(t)=W_m^{-1/2}\mathbf z(t).
\end{equation*}
The eigenvalues of $S_m$ determine the decay rates of the weighted error along the corresponding eigenvectors. These eigenvectors generally combine the two Fourier modes, so each Fourier-coefficient error is a linear combination of the corresponding exponential terms.

To compare the decay rates of the individual Fourier modes, we additionally assume that $\overline G$ is diagonal in the prescribed Fourier basis,
\begin{equation}
  \overline G=\operatorname{diag}\bigl(\lambda_{\mathrm{net}}(1),\lambda_{\mathrm{net}}(k)\bigr),
  \qquad
  \lambda_{\mathrm{net}}(1)>0,\quad \lambda_{\mathrm{net}}(k)\geq0.
  \label{eq:appendix-diagonal-kernel}
\end{equation}
This assumption eliminates coupling between the Fourier coefficients and is not implied by the spectral decomposition above. The errors then satisfy $e_j(t)=e_j(0)\exp(-r_m(j)t)$ with $r_m(j)=\lambda_{\mathrm{net}}(j)(j\pi)^{2m}$ for $j\in\{1,k\}$, and hence
\begin{equation}
  \frac{r_m(k)}{r_m(1)}
  =\frac{\lambda_{\mathrm{net}}(k)}{\lambda_{\mathrm{net}}(1)}k^{2m}.
  \label{eq:appendix-rate-ratio}
\end{equation}
Relative to regression, the Deep Ritz and residual objectives multiply this ratio by $k^2$ and $k^4$, respectively. If $\lambda_{\mathrm{net}}(k)\asymp k^{-p}$ as $k\to\infty$, with $p>0$, then $r_m(k)\asymp k^{2m-p}$. For a fixed positive low-frequency rate, the high-to-low rate ratio therefore tends to zero when $2m<p$, remains of order one when $2m=p$, and diverges when $2m>p$. Thus, in this uncoupled model, differential-loss weighting can offset or reverse the asymptotic decay of the modal rates with frequency.

\clearpage
\section{Image regression training curves}
\label{app:image-regression-curves}
\begin{figure}[H]
  \centering
  \includegraphics[width=\linewidth]{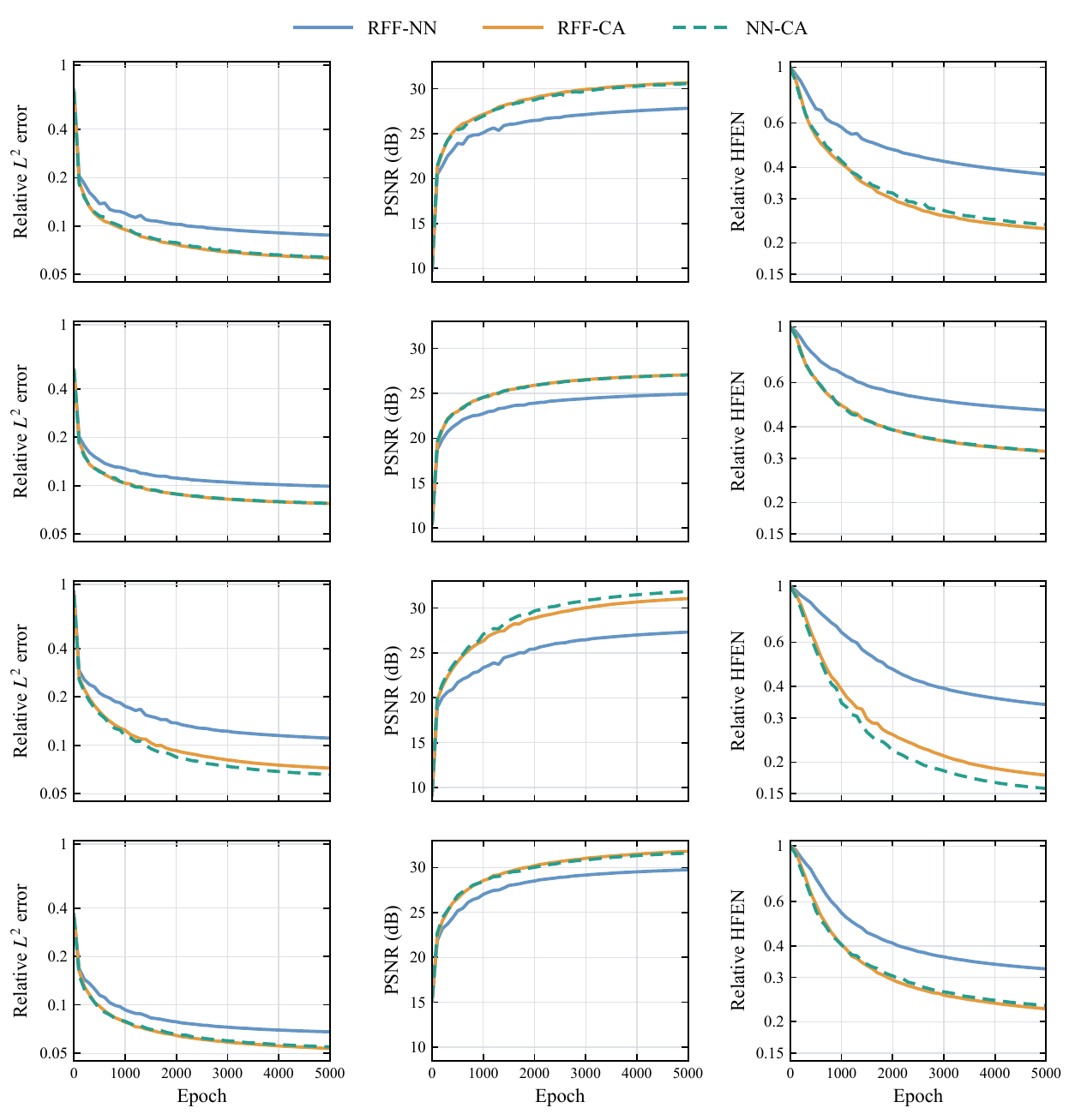}
  \caption{Relative $L^2$ error (left), PSNR (middle), and relative HFEN (right) during training for RFF-NN, RFF-CA, and NN-CA on four DIV2K images. From top to bottom, rows correspond to the images shown from left to right in Fig.~\ref{fig:div2k-samples}.}
  \label{fig:img_results}
\end{figure}

Figure~\ref{fig:img_results} shows the training curves for the four DIV2K image experiments whose final results are reported in Table~\ref{tab:img_results}.

\clearpage
\bibliographystyle{unsrt}
\bibliography{ref}
\end{document}